\documentclass[a4paper,10pt]{article}
\usepackage{color}   
\usepackage[bookmarks=true,pdfencoding=unicode]{hyperref}
\hypersetup{
    backref=true,
    citecolor=blue,
    colorlinks=true, 
    linktoc=all,     
    linkcolor=blue,  
}

\usepackage{bookmark}
\usepackage[latin1]{inputenc}
\usepackage[T1]{fontenc}
\usepackage{bezier}
\usepackage{dsfont}
\usepackage{mathptmx}      
\usepackage{amsfonts}
\usepackage{amsmath}
\usepackage{amssymb}
\usepackage{amsbsy}
\usepackage{amscd}
\usepackage{amsopn}
\usepackage{amssymb}
\usepackage{amstext}

\usepackage{amsxtra}
\usepackage{array}
\usepackage{bm}
\usepackage{graphicx,psfrag}
\usepackage{latexsym}
\usepackage{multirow}
\usepackage{subfig}
\usepackage{etoolbox}
\usepackage{algorithmicx}
\usepackage{algpseudocode}
\algtext*{EndIf}
\usepackage{mdframed}
\usepackage{lipsum} 

\newcommand{\wqed}{\hfill \ensuremath{\Box}}

\DeclareMathAlphabet{\mathcal}{OMS}{cmsy}{m}{n}

\newtoggle{LatexWithProofs}
\togglefalse{LatexWithProofs}

\newtoggle{LatexFull}
\togglefalse{LatexFull}

\newtheorem{pLemma}{Lemma}
\newtheorem{pProp}{Proposition}
\newtheorem{pArgument}{Argument}
\newtheorem{pExample}{Example}
\newtheorem{pCorol}{Corollary}
\newtheorem{pLongDef}{Definition}
\newtheorem{pTheorem}{Theorem}
\newtheorem{pConjecture}{Conjecture}

\newcommand{\pRef}[1]{(\ref{#1})}

\newcommand{\pbFigureBH}[2]{
\begin{figure}[h]
\label{#1}
\hypertarget{#1}{}
\bookmark[
rellevel=1,
keeplevel,
dest=#1
]{Fig \ref{#1}: {#2}}
}

\newcommand{\peFigure}[1]{
\label{#1}
\end{figure}}

\newcommand{\peLongDef}[1]{$\blacktriangle$ \end{pLongDef}}
\newcommand{\peLongDefX}[1]{\end{pLongDef}}

\newcommand{\peArgument}[1]{ \hyperlink{#1}{$\blacktriangle$} \end{pArgument}}

\newcommand{\peExample}[1]{ \hyperlink{#1}{$\blacktriangle$} \end{pExample}}
\newcommand{\peExampleX}[1]{\end{pExample}}

\newcommand{\pePlainExample}[1]{$\blacktriangle$ \end{pExample}}
\newcommand{\pePlainExampleX}[1]{\end{pExample}}

\newcommand{\peCorol}[1]{\hyperlink{#1}{$\blacktriangle$} \end{pCorol}}
\newcommand{\peCorolX}[1]{\end{pCorol}}

\newcommand{\peLemma}[1]{ \hyperlink{#1}{$\blacktriangle$} \end{pLemma}}
\newcommand{\peLemmaX}[1]{\end{pLemma}}

\newcommand{\peProp}[1]{ \hyperlink{#1}{$\blacktriangle$} \end{pProp}}
\newcommand{\pePropX}[1]{\end{pProp}}

\newcommand{\peTheorem}[1]{ \hyperlink{#1}{$\blacktriangle$} \end{pTheorem}}
\newcommand{\peConjecture}[1]{ \hyperlink{#1}{$\blacktriangle$} \end{pConjecture}}

\iftoggle{LatexFull}
{
\newcommand{\pFullLink}[1]{\eqno \hyperlink{#1}{\blacktriangle}}
\newcommand{\peFullProp}[1]{ \hyperlink{#1}{$\blacktriangle$} \end{pProp}}
\newcommand{\peFullCorol}[1]{ \hyperlink{#1}{$\blacktriangle$} \end{pCorol}}
\newcommand{\peFullLemma}[1]{ \hyperlink{#1}{$\blacktriangle$} \end{pLemma}}
\newcommand{\peFullExample}[1]{ \hyperlink{#1}{$\blacktriangle$} \end{pExample}}
}
{
\newcommand{\pFullLink}[1]{\eqno \blacktriangle}
\newcommand{\peFullProp}[1]{$\blacktriangle$ \end{pProp}}
\newcommand{\peFullCorol}[1]{$\blacktriangle$ \end{pCorol}}
\newcommand{\peFullLemma}[1]{$\blacktriangle$ \end{pLemma}}
\newcommand{\peFullExample}[1]{$\blacktriangle$ \end{pExample}}
}

\iftoggle{LatexWithProofs}
{
\definecolor{darkpink}{rgb}{0.91, 0.33, 0.5}
\definecolor{darksalmon}{rgb}{0.91, 0.59, 0.48}
\definecolor{desertsand}{rgb}{0.93, 0.79, 0.69}
\definecolor{celadon}{rgb}{0.67, 0.88, 0.69}
\definecolor{darkcyan}{rgb}{0.0, 0.55, 0.55}

\newcommand{\pbClaim}[1]{{\color{red} ?}}
\newcommand{\peClaim}[1]{{\color{red} ?`}}
\newcommand{\pbArgument}[1]{\begin{pArgument} \label{#1}  {\color{darksalmon} #1}}

\newcommand{\pbExample}[1]{\begin{pExample} \label{#1}  {\color{darksalmon} #1}}

\newcommand{\pbExampleB}[1]{
\begin{pExample} \label{#1}
{\color{darksalmon} #1}
\hypertarget{{T#1}}{}
\bookmark[
rellevel=1,
keeplevel,
dest=T#1
]{Example \ref{#1}}
}

\newcommand{\pbExampleBT}[2]{
\begin{pExample}[#2] \label{#1}
{\color{darksalmon} #1}
\hypertarget{{T#1}}{}
\bookmark[
rellevel=1,
keeplevel,
dest=T#1
]{Example \ref{#1}: {#2}}
}

\newcommand{\pbCorol}[1]{\begin{pCorol} \label{#1}  {\color{darksalmon} #1}}

\newcommand{\pbCorolB}[1]{
\begin{pCorol} \label{#1}  {\color{darksalmon} #1}
\hypertarget{{T#1}}{}
\bookmark[
rellevel=1,
keeplevel,
dest=T#1
]{Corollary \ref{#1}}
}

\newcommand{\pbCorolBT}[2]{
\begin{pCorol}[#2]\label{#1}  {\color{darksalmon} #1}
\hypertarget{{T#1}}{}
\bookmark[
rellevel=1,
keeplevel,
dest=T#1
]{Corollary \ref{#1}: {#2}}
}

\newcommand{\pbLemma}[1]{\begin{pLemma} \label{#1}  {\color{darksalmon} #1}}

\newcommand{\pbLemmaB}[1]{
\begin{pLemma} \label{#1}  {\color{darksalmon} #1}
\hypertarget{{T#1}}{}
\bookmark[
rellevel=1,
keeplevel,
dest=T#1
]{Lemma \ref{#1}}
}

\newcommand{\pbLemmaBT}[2]{
\begin{pLemma}[#2] \label{#1}  {\color{darksalmon} #1}
\hypertarget{{T#1}}{}
\bookmark[
rellevel=1,
keeplevel,
dest=T#1
]{Lemma \ref{#1}: {#2}}
}

\newcommand{\pbProp}[1]{\begin{pProp} \label{#1}  {\color{darksalmon} #1}}

\newcommand{\pbPropB}[1]{
\begin{pProp} \label{#1}  {\color{darksalmon} #1}
\hypertarget{{T#1}}{}
\bookmark[
rellevel=1,
keeplevel,
dest=T#1
]{Prop. \ref{#1}}
}

\newcommand{\pbPropBT}[2]{
\begin{pProp}[#2] \label{#1}  {\color{darksalmon} #1}
\hypertarget{{T#1}}{}
\bookmark[
rellevel=1,
keeplevel,
dest=T#1
]{Prop. \ref{#1}: {#2}}
}

\newcommand{\pbTheorem}[1]{\begin{pTheorem} \label{#1}  {\color{darksalmon} #1}}

\newcommand{\pbTheoremB}[1]{
\begin{pTheorem} \label{#1}  {\color{darksalmon} #1}
\hypertarget{{T#1}}{}
\bookmark[
rellevel=1,
keeplevel,
dest=T#1
]{Theorem \ref{#1}}
}

\newcommand{\pbTheoremBT}[2]{
\begin{pTheorem}[#2] \label{#1}  {\color{darksalmon} #1}
\hypertarget{{T#1}}{}
\bookmark[
rellevel=1,
keeplevel,
dest=T#1
]{Theorem \ref{#1}: {#2}}
}

\newcommand{\pbConjectureB}[1]{
\begin{pConjecture} \label{#1}  {\color{darksalmon} #1}
\hypertarget{{T#1}}{}
\bookmark[
rellevel=1,
keeplevel,
dest=T#1
]{Conjecture \ref{#1}}
}

\newcommand{\pbConjectureBT}[2]{
\begin{pConjecture}[#2] \label{#1}  {\color{darksalmon} #1}
\hypertarget{{T#1}}{}
\bookmark[
rellevel=1,
keeplevel,
dest=T#1
]{Conjecture \ref{#1}: {#2}}
}

\newcommand{\pbLongDef}[1]{\begin{pLongDef}
\label{#1} \hypertarget{{#1}}{} {\color{darksalmon} #1}
\bookmark[
rellevel=1,
keeplevel,
dest=#1
]{Definition \ref{#1}}
}

\newcommand{\pbLongDefB}[2]{\begin{pLongDef}
\label{#1} \hypertarget{{#1}}{} {\color{darksalmon} #1}
\bookmark[
rellevel=1,
keeplevel,
dest=#1
]{\ref{#1}: #2}
}

\newcommand{\pbLongDefBT}[2]{\begin{pLongDef}[#2]
\label{#1} \hypertarget{{#1}}{} {\color{darksalmon} #1}
\bookmark[
rellevel=1,
keeplevel,
dest=#1
]{\ref{#1}: #2}
}

\newcommand{\pbChain}[1]{\[ }
\newcommand{\peChain}[1]{\ {\color{red} \checkmark  \wrm{#1}} \] }

\newcommand{\pePClaim}[1]{ \checkmark}

\newcommand{\pbTClaim}[1]{\begin{equation} \label{#1}  }
\newcommand{\peTClaim}[1]{ \ {\color{red} \checkmark  \wrm{#1}} \end{equation} }
\newcommand{\peEClaim}[1]{ \ {\color{red} \checkmark  \wrm{#1 \bullet} }\end{equation} }

\newcommand{\pbDef}[1] {\hypertarget{{#1}}{} \begin{equation} \label{#1} }
\newcommand{\peDef}[1]{ {\color{blue} \bigstar \wrm{#1}} \end{equation}}
\newcommand{\pDefLabel}[1]{\label{#1} {\color{blue} \bigstar \wrm{#1}}}

\newcommand{\pbHypot}[1]{ \begin{equation} \label{#1}  }
\newcommand{\peHypot}[1]{  \ \ {{\color{blue} \bigstar \wrm{#1}}}  \end{equation} }

\newcommand{\pbProofB}[2]{\newpage {\color{darksalmon} Proof of {#1} {#2}.} \ref{#2} \hypertarget{{#2}}{}
\bookmark[
rellevel=1,
keeplevel,
dest=#2
]
{Proof of {#1} \ref{#2}}
}

\newcommand{\peProof}[2]{\wqed{} {\color{darksalmon} End of proof of #1 #2} \ref{#2}}

\newcommand{\peVerify}[1]{\wqed{} {\color{darksalmon} End of verification of Example #1} \ref{#1}}
}
{
\newcommand{\pbChain}[1]{\[}
\newcommand{\peChain}[1]{\]}

\newcommand{\pbClaim}[1]{}
\newcommand{\peClaim}[1]{}

\newcommand{\peTClaim}[1]{ \end{equation}}
\newcommand{\peEClaim}[1]{ \end{equation}}
\newcommand{\pbTClaim}[1]{\begin{equation} \label{#1}}

\newcommand{\peDef}[1]{  \end{equation}}

\newcommand{\lpPlainClaim}[1]{}
\newcommand{\lpEndPlainClaim}[1]{}
\newcommand{\pbDef}[1]{ \hypertarget{{#1}}{} \begin{equation} \label{#1}}
\newcommand{\pDefLabel}[1]{\label{#1}}
\newcommand{\pbHypot}[1]{\begin{equation} \label{#1}  }
\newcommand{\peHypot}[1]{ \end{equation}}

\newcommand{\pbCorol}[1]{\begin{pCorol} \label{#1}}

\newcommand{\pbCorolB}[1]{\begin{pCorol} \label{#1}
\hypertarget{{T#1}}{}
\bookmark[
rellevel=1,
keeplevel,
dest=T#1
]
{Corollary \ref{#1}}
}

\newcommand{\pbCorolBT}[2]{
\begin{pCorol}[#2] \label{#1}
\hypertarget{{T#1}}{}
\bookmark[
rellevel=1,
keeplevel,
dest=T#1
]
{Corollary \ref{#1}: {#2}}
}

\newcommand{\pbArgument}[1]{\begin{pArgument} \label{#1}}
\newcommand{\pbExample}[1]{\begin{pExample} \label{#1}}

\newcommand{\pbExampleB}[1]{\begin{pExample} \label{#1}
\hypertarget{{T#1}}{}
\bookmark[
rellevel=1,
keeplevel,
dest=T#1
]{Example \ref{#1}}
}

\newcommand{\pbExampleBT}[2]{
\begin{pExample}[#2] \label{#1}
\hypertarget{{T#1}}{}
\bookmark[
rellevel=1,
keeplevel,
dest=T#1
]{Example \ref{#1}: {#2}}
}

\newcommand{\pbLemma}[1]{\begin{pLemma} \label{#1}}

\newcommand{\pbLemmaBT}[2]{\begin{pLemma}[#2] \label{#1}
\hypertarget{{T#1}}{}
\bookmark[
rellevel=1,
keeplevel,
dest=T#1
]
{Lemma \ref{#1}: {#2}}
}

\newcommand{\pbConjectureBT}[2]{
\begin{pConjecture}[#2] \label{#1}
\hypertarget{{T#1}}{}
\bookmark[
rellevel=1,
keeplevel,
dest=T#1s
]{Conjecture \ref{#1}: {#2}}
}

\newcommand{\pbTheoremBT}[2]{
\begin{pTheorem}[#2] \label{#1}
\hypertarget{{T#1}}{}
\bookmark[
rellevel=1,
keeplevel,
dest=T#1
]{Theorem \ref{#1}: {#2}}
}

\newcommand{\pbLemmaB}[1]{\begin{pLemma} \label{#1}
\hypertarget{{T#1}}{}
\bookmark[
rellevel=1,
keeplevel,
dest=T#1
]
{Lemma \ref{#1}}
}

\newcommand{\pbProp}[1]{\begin{pProp} \label{#1}}
\newcommand{\pbPropBT}[2]{\begin{pProp}[#2] \label{#1}
\hypertarget{{T#1}}{}
\bookmark[
rellevel=1,
keeplevel,
dest=T#1
]
{Prop. \ref{#1}: {#2}}
}

\newcommand{\pbPropB}[1]{\begin{pProp} \label{#1}
\hypertarget{{T#1}}{}
\bookmark[
rellevel=1,
keeplevel,
dest=T#1
]
{Prop. \ref{#1}}
}

\newcommand{\pbTheorem}[1]{\begin{pTheorem} \label{#1}}
\newcommand{\pbLongDef}[1]{\begin{pLongDef} \label{#1}}

\newcommand{\pbLongDefB}[2]{\begin{pLongDef}
\label{#1} \hypertarget{{#1}}{}
\bookmark[
rellevel=1,
keeplevel,
dest=#1
]{\ref{#1}: #2}
}

\newcommand{\pbLongDefBT}[2]{\begin{pLongDef}[#2]
\label{#1} \hypertarget{{#1}}{}
\bookmark[
rellevel=1,
keeplevel,
dest=#1
]{\ref{#1}: #2}
}

\newcommand{\pbProofB}[2]{ {\bf Proof of {#1} \ref{#2}.} \hypertarget{{#2}}{}
\bookmark[
rellevel=1,
keeplevel,
dest=#2
]
{Proof of  {#1} \ref{#2}}
}

\newcommand{\peProof}[2]{\wqed{}}

\newcommand{\peVerify}[1]{\wqed{}}

}

\newcommand{\wabs}[1]{\left|#1\right|}

\newcommand{\wc}[1]{\mathds{C}}
\newcommand{\wcal}[1]{\mathcal{#1}}

\newcommand{\wceilB}[1]{\Big{\lceil} #1 \Big{\rceil}}

\newcommand{\wdfc}[2]{{#1}'\!\left(#2\right)}

\newcommand{\wdsf}[2]{{#1}''\!\!\left(#2\right)}

\newcommand{\wfc}[2]{{#1}\!\left(#2\right)}

\newcommand{\wfloor}[1]{\lfloor {{#1}} \rfloor }

\newcommand{\wfloorB}[1]{\Big{\lfloor} {{#1}} \Big{\rfloor}}

\newcommand{\wi}[1]{\wrm{i}}

\algdef{SE}[DOWHILE]{Do}{doWhile}{\algorithmicdo}[1]{\algorithmicwhile\ #1}%

\newdimen\CdotAxis
\newcommand*{\CdotAux}[3]{%
  {%
    \settoheight\CdotAxis{$#2\vcenter{}$}%
    \sbox0{%
      \raisebox\CdotAxis{%
        \scalebox{#1}{%
          \raisebox{-\CdotAxis}{%
            $\mathsurround=0pt #2#3$%
          }%
        }%
      }%
    }%
    \dp0=0pt %
    \sbox2{$#2\bullet$}%
    \ifdim\ht2<\ht0 %
      \ht0=\ht2 %
    \fi
    \sbox2{$\mathsurround=0pt #2#3$}%
    \hbox to \wd2{\hss\usebox{0}\hss}%
  }%
}

\newcommand{\wlr}[1]{\left( #1 \right)}

\newcommand{\wn}{\mathds N}

\newcommand{\wnorm}[1]{\left\| #1 \right\|}

\newcommand{\wrone}{\mathds R}

\newcommand{\wrm}[1]{\mathrm{#1}}

\newcommand{\wset}[1]{{\left\{ #1 \right\}}}

\newcommand{\wz}{\mathds Z}

\newcommand{\wfpes}[1]{\wcal{E}_{\wcal{A},e}}

\newcolumntype{M}[1]{>{\centering\arraybackslash}m{#1}}
\newcolumntype{N}{@{}m{0pt}@{}}
\begin{document}
\title{On the rate of convergence of Berrut's interpolant at equally spaced nodes}

\author{Walter F. Mascarenhas\thanks{
Instituto de Matem\'{a}tica e Estat\'{i}stica, Universidade de S\~{a}o Paulo,
Cidade Universit\'{a}ria, Rua do Mat\~{a}o 1010, S\~{a}o Paulo SP, Brazil. CEP 05508-090.
Tel.: +55-11-3091 5411, Fax: +55-11-3091 6134, walter.mascarenhas@gmail.com.}}
\markboth{W.F. Mascarenhas}{Stability versus order of approximation in barycentric interpolation}
\maketitle

\begin{abstract}
We extend the  work by Mastroianni and Szabados regarding
the barycentric interpolant introduced by J.-P. Berrut in 1988,
for equally spaced nodes. We prove fully their first conjecture
and present a proof of a weaker version of their second conjecture.
More importantly than proving these conjectures,  
we present a sharp description of the asymptotic error incurred by
the interpolants when the derivative of the interpolated function is 
absolutely continuous, which is a class of functions 
broad enough to cover most functions usually found in practice. 
We also contribute to the solution 
of the broad problem they raised regarding the order of approximation of these interpolants, 
by showing that they have
order of approximation of order 1/n for functions
with derivatives of bounded variation.
\end{abstract}

\section{Introduction}
\label{secIntro}
In a recent article \cite{MastroSz},
professors G. Mastroianni and J. Szabados
discuss barycentric interpolation
of functions $f: [-1,1] \rightarrow \wrone{}$
at equally spaced nodes
\[
x_{k,n} := 2 k / n - 1 \ \wrm{for} \ \ k = 0,\dots,n.
\]
They analyze the order of approximation of the
barycentric interpolant
introduced by J.-P. Berrut \cite{BerrutRat}:
\pbDef{bary}
\wfc{B_{n}}{f,x} := \frac{\wfc{N_n}{f,x}}{\wfc{D_n}{x}}
\hspace{0.3cm}  \wrm{for} \hspace{0.3cm}
\ x \not \in \wset{x_{0,n},\dots,x_{n,n}}
\hspace{0.3cm}
\wrm{and}
\hspace{0.3cm}
\wfc{B_n}{f,x_{k,n}} = \wfc{f}{x_{k,n}},
\peDef{bary}
with
\pbDef{numDen}
\wfc{N_n}{f,x} := \sum_{k = 0}^n \wlr{-1}^k \frac{\wfc{f}{x_{k,n}}}{x - x_{k,n}}
\hspace{1cm} \wrm{and} \hspace{1cm}
\wfc{D_n}{x} := \sum_{k = 0}^n \wlr{-1}^k \frac{1}{x - x_{k,n}}.
\peDef{numDen}
They proved some results and stated two conjectures
and a broad open problem about the rate at which $\wfc{B_n}{f}$
approximates $f$ for some classes of functions.
Although their proof of their second theorem is incorrect,
their conclusions are correct and they have correctly shown 
that the error
$\wnorm{\wfc{B_n}{f} - f}_\infty$ is of order
$1/n$ for functions with derivatives in
the class $\wrm{Lip \ 1}$ of
functions with continuity modulus
$\wfc{\omega}{t} \leq \kappa t$.

In this article we extend their work,
by presenting a detailed analysis of the asymptotic
behavior of the interpolation error for functions
with absolutely continuous derivatives.
We denote the class of such functions by $\wrm{AC}^1$,
and emphasize that, unlike
the definition of the Sobolev space $\wfc{\wrm{W}^{2,1}}{[-1,1]}$,
we require that $\wdfc{f}{x}$ is defined
for all $x \in [-1,1]$ in order for $f$ to
belong to $\wrm{AC}^1$ (we consider directional derivatives
at $x \in \wset{-1,1}$.) We also analyze functions with
derivatives of bounded variation, and denote their class
by $\wrm{BV}^1$, with the same requirement on the derivatives.

We prove the first conjecture by Mastroianni and Szabados
in full, and present a proof of a weaker version of their
second conjecture: their conjecture regards arbitrary functions,
our proof makes the additional assumption that the function
have absolutely continuous derivatives, but we hope that
the readers will agree with us that this class of functions
covers a wide range of applications. We also show that the order 
of convergence of the interpolants $\wfc{B_n}{f}$ above
is also of order $1/n$ for $f \in \wrm{BV}^1$. 
Their first conjecture, which we state below,
is about the interpolation error for
functions $f \in \wrm{Lip \ 1}$,
and we prove it in Section \ref{secFirst}.
\pbConjectureBT{conj}{First conjecture by Mastroianni and Szabados}
There exists a function $f \in \wrm{Lip \ 1}$ such that
\[
\limsup_{n \rightarrow \infty} \frac{n}{\log n} \wnorm{\wfc{B_{n}}{f}- f}_\infty > 0.
\]
\peConjecture{conj}

Regarding the second conjecture, we have found that if 
$f \in \wrm{BV}^1$  then we can bound the sequence
\[
n \ \wnorm{\wfc{B}{f} - f}_{\infty}
\]
by a constant depending on $f$. Moreover, if
$f \in \wrm{AC}^1$ and $x \in [-1,1]$ then we
can describe exactly all possible accumulation points
of the sequence
\pbDef{error}
n \ \wlr{\wfc{B_n}{f,x} - \wfc{f}{x}}.
\peDef{error}

This description is given by
Theorem \ref{thmMain} below and uses the functions
\pbDef{oddLimit}
\wfc{O}{f,x} :=
\frac{\wfc{f}{x} - \wfc{f}{1}}{2 \wlr{x - 1}} - \frac{\wfc{f}{x} - \wfc{f}{-1}}{2 \wlr{x + 1}},
\peDef{oddLimit}
\pbDef{evenLimit}
\wfc{E}{f,x} :=
\frac{\wfc{f}{1} - \wfc{f}{x}}{2\wlr{x - 1}} + \frac{\wfc{f}{-1} - \wfc{f}{x}}{2 \wlr{x + 1}}.
\peDef{evenLimit}
(Throughout the article,
$O$ stands for {\it odd} and $E$ stands for {\it even}.)

We must be careful when analyzing the
sequences in Equation \pRef{error} when $x = x_{k,n}$ is a node,
because both the denominator and the numerator
of $\wfc{B_n}{f,x}$
are discontinuous at such $x$, and the interpolant
is defined in a different way for them
in Equation \pRef{bary}.
As a result, the error has more favourable properties
at the nodes and this may confuse our analysis of
the convergence for a general $x$.
For instance,  if $x \in \wset{-1,1}$ then the error
$\wfc{B_n}{f,x} - \wfc{f}{x}$ is zero for all $n$,
and the same holds for $x = 0$ when $n$ is even.
In order to handle this issue precisely,
 we state the following definitions:

\pbLongDefBT{defSequence}{Sequence}
We say that an increasing function
$n: \wn{} \rightarrow \wn{}$ with $\wfc{n}{j} = n_j$
is ``a sequence $n_j$.'' The sequence is
odd if $n_j$ is odd for all $j$, and it
is even if $n_j$ is even for all $j$.
\peLongDef{defRegular}

\pbLongDefBT{defRegular}{Regular point}
We say that $x \in [-1,1]$ is
regular for the sequence $n_j$ if there exists $j_0$ such
that
\[
j \geq j_0 \ \Rightarrow \ x \not \in \wset{x_{0,n_j}, \, x_{1,n_j}, \, \dots, \, x_{n_j,n_j} \, }.
\]
\peLongDef{defRegular}

\pbLongDefBT{defRext}{The compactification of $\wrone{}$}
In order to handle infinite limits, we
write
\[
\overline{\wrone{}} := \wrone{} \ \bigcup \ \wset{+\infty,-\infty}
\]
as the two point compactification of $\wrone{}$, endowed
with the usual topology and extension of the operators $<$ and $\leq$.
In particular, $\overline{\wrone{}}$ and its subset $[\pi/2,+\infty]$,
which are relevant to our discussion,
are compact in our topology for $\overline{\wrone{}}$.
\peLongDef{defRExt}

All irrational points are regular for every sequence $n_j$;
the points $\pm 1$ are not regular for any sequence, and
$0$ is regular for odd sequences and irregular for even ones.
Given $x \in [-1,1]$, we can decompose any sequence
$n_j$ in at most three parts: one in which
$x$ is a node for all $j$, so that
$\wfc{B_{n_j}}{f,x} = \wfc{f}{x}$ for all $j$,
an two other sequences for which $x$ is regular,
one even and another odd
(of course, some parts
may not be necessary.) Therefore, by understanding
the regular points for odd and even sequences we
can get the full picture regarding the pointwise
convergence of the interpolation error.
We now state our first formal result.

%

\pbTheoremBT{thmMain}{The limits of $n \, \wlr{\wfc{B_n}{f,x} - \wfc{f}{x}}$ for $f$ in $\wrm{AC}^1$}
Let $f$ be a function in $\wrm{AC}^1$,
$n_j$ an odd sequence,  and $x \in [-1,1]$
such that
\pbDef{converge}
\lim_{j \rightarrow \infty} \ n_j \ \wlr{\wfc{B_{n_j}}{f,x} - \wfc{f}{x}} = L \in \overline{\wrone{}}.
\peDef{converge}
If $x$ is irrational then, for the function $\wfc{O}{f,x}$ in Equation \pRef{oddLimit},
\pbTClaim{ofIrrat}
L \in \wfc{\wcal{O}}{f,x} := \left[ -\frac{2 \wabs{\wfc{O}{f,x}}}{\pi},
                                   \ \frac{2 \wabs{\wfc{O}{f,x}}}{\pi} \right],
\peTClaim{ofIrrat}
and if $x$ is rational then there exists a finite set
$\wfc{\wcal{O}}{x} \subset \wrone{} \setminus \wset{0}$,
defined in Equation \pRef{oddSetRat} in Section \ref{secSecond},
such that and if $x$  is a regular rational point for
$n_j$ then
\[
L \in \wfc{\wcal{O}}{f,x}
:= \wset{  \wfc{O}{f,x} / y, \ y \in \wfc{\wcal{O}}{x} }.
\]
Conversely, if $L \in \wfc{\wcal{O}}{f,x}$ then
there exists an odd sequence $n_j$ for which
$x$ is regular and Equation \pRef{converge} holds.

Similarly, if $n_j$ is an even sequence,
Equation \pRef{converge} holds and
$x$ is irrational then
\[
L \in \wfc{\wcal{E}}{f,x} := \left[ -\frac{2 \wabs{\wfc{E}{f,x}}}{\pi},
                                   \ \frac{2 \wabs{\wfc{E}{f,x}}}{\pi} \right],
\]
and if $x$ is rational then there exists a finite set
$\wfc{\wcal{E}}{x} \subset \wrone{} \setminus \wset{0}$,
defined in Equation \pRef{evenSetRat} in Section \ref{secSecond},
such that and if $x$  is a regular rational point for
$n_j$ then
\[
L \in \wfc{\wcal{E}}{f,x}
:= \wset{\wfc{E}{f,x}  / y, \ y \in \wfc{\wcal{E}}{x} }.
\]
Conversely,
if $L \in \wfc{\wcal{E}}{f,x}$ then
there exists an even sequence $n_j$
for which $x$ is a regular point and Equation \pRef{converge} holds.
\peTheorem{thmMain}

Theorem \ref{thmMain} has far reaching implications for $f \in \wrm{AC}^1$.
For instance, it yields a simple proof of
second conjecture by Mastroianni and Szabados stated below, 
with the additional hypothesis that $f$ is in this class:
\pbConjectureBT{conj2}{Second conjecture by Mastroianni and Szabados}
We have
\[
\wnorm{\wfc{B_n}{f} - f}_\infty = \wfc{o}{1/n}
\]
if and only if $f$ is constant
(when $n =2,4,...$), or $f$ is linear (when $n = 1,3...$).
\peConjecture{conj2}

In fact, when $f \in \wrm{AC}^1$,
if $\wnorm{\wfc{B_{n}}{f} - f}_\infty = \wfc{o}{1/n}$
and $z \in (-1,1)$ is irrational then  Theorem \ref{thmMain} implies that
$\wfc{\wcal{O}}{f,z} = \wset{0}$ and Equation \pRef{oddLimit} leads to
\pbTClaim{irrat}
\wfc{O}{f,z} =
\frac{\wfc{f}{z} - \wfc{f}{1}}{2 \wlr{z - 1}} - \frac{\wfc{f}{z} - \wfc{f}{-1}}{2 \wlr{z + 1}} = 0,
\peTClaim{irrat}
and by the continuity of $f$
Equation \pRef{irrat} must hold for all $x \in [-1,1]$. Therefore,
\[
\wfc{f}{x} = \frac{\wfc{f}{1} + \wfc{f}{-1}}{2} + \frac{\wfc{f}{1} - \wfc{f}{-1}}{2} x,
\]
and $f$ is linear. This proves the second conjecture for odd sequences.

The same argument using the part of Theorem \ref{thmMain} for even sequences
leads to
\[
\wfc{E}{f,x} = \frac{\wfc{f}{1} - \wfc{f}{x}}{2 \wlr{x - 1}} + \frac{\wfc{f}{-1} - \wfc{f}{x}}{2 \wlr{x + 1}} = 0.
\]
For $x \neq 0$ this equation implies that
\pbDef{constf}
\wfc{f}{x} = \frac{\wlr{x - 1} \wfc{f}{-1} + \wlr{x + 1}\wfc{f}{1}}{2 x},
\peDef{constf}
the continuity of $f$ at $x = 0$ yields $\wfc{f}{1} = \wfc{f}{-1}$,
and Equation \pRef{constf} shows that $f$ is constant. This
finishes the proof of the second conjecture for $f \in \wrm{AC}^1$.

Besides the weakened version of the second conjecture above,
we can prove other interesting results using
Theorem \ref{thmMain}. For instance,
if $x$ is rational then
$0 \not \in \wfc{\wcal{O}}{x} \bigcup \wfc{\wcal{E}}{x}$
and the reader will be able to prove the following corollary:

\pbCorolBT{corBarRat}{Large errors for rational $x$}
If $f \in \wrm{AC}^1$ and
$x \in [-1,1]$ is rational and regular for the sequence $n_j$,
$\wfc{O}{f,x} \neq 0$
and $\wfc{E}{f,x} \neq 0$ then
\[
\liminf_{j \rightarrow \infty} \ n_j \ \wabs{\wfc{B_{n_j}}{f,x} - \wfc{f}{x}} > 0.
\]
\peCorol{corBarRat}

However, Theorem \ref{thmMain} has a serious limitation: it is
only a pointwise result, and it does not imply the more
interesting bound
\[
\limsup_{n \rightarrow \infty} \ n \, \wnorm{\wfc{B_n}{f} - f}_\infty  < +\infty
\]
considered by Mastroianni and Szabados in their open problem.
Fortunately, we can also prove uniform convergence 
results for $f \in \wrm{AC}^1$:

\pbTheoremBT{thmUnif}{Uniform convergence for $f \in \wrm{AC}^1$}
If $f \in \wrm{AC}^1$ and $n_j$ is an odd sequence
then, for the function $\wfc{O}{f}$ defined in Equation \pRef{oddLimit},
\[
\lim_{j \rightarrow \infty}  \  n_j  \ \wnorm{\wfc{B_{n_j}}{f} - f - \wfc{O}{f}/D_{n_j}}_\infty = 0
\]
and if $n_j$ is an even sequence then
\[
\lim_{j \rightarrow \infty} \  n_j \ \wnorm{\wfc{B_{n_j}}{f} - f - \wfc{E}{f}/D_{n_j}}_\infty = 0,
\]
for $\wfc{E}{f}$ defined in Equation \pRef{evenLimit},
\peTheorem{thmUnif}
Lemma \ref{lemDen} in Section \ref{secDen}
yields $\wnorm{n/D_n}_\infty \leq 1$, and it is
clear that $\wnorm{\wfc{O}{f}}_\infty \leq \wnorm{f'}_\infty$
and $\wnorm{\wfc{E}{f}}_\infty \leq \wnorm{f'}_\infty$ .
These observations combined with
Theorem \ref{thmUnif} lead to an uniform upper bound of order $1/n$ in the
interpolation error for $f \in \wrm{AC}^1$, but we can derive
this bound under the weaker assumption of derivatives of bounded variation:
\pbTheoremBT{thmUnifBV}{Uniform convergence when $f \in \wrm{BV}^1$}
If $f \in \wrm{BV}^1$ then 
\pbDef{unifBV}
n \wnorm{\wfc{B_{n}}{f} - f}_\infty \leq 
T_{f'}[-1,1] /2 + \max \wset{\wnorm{\wfc{O}{f}}_\infty, \wnorm{\wfc{E}{f}}_\infty},
\peDef{unifBV}
where $T_{f'}[-1,1]$ is the total variation of $f'$ in $[-1,1]$.
\peTheorem{thmUnifBV}

We prove the results above in the next sections.
In Section \ref{secFirst} we prove the first conjecture.
In Section \ref{secDen} we discuss the denominator
of the interpolant defined in Equation \pRef{bary}.
In Section \ref{secNum} we analyze
the numerator of the error $\wfc{B_n}{f,x} - \wfc{f}{x}$
for functions in $\wrm{AC}^1$. In Section \ref{secBV}
we analyze the numerator for $f \in \wrm{BV}^1$.
Finally, in Section \ref{secSecond} we combine the results
in Sections \ref{secDen}. \ref{secNum} and \ref{secBV} to
prove Theorems \ref{thmMain}, 
\ref{thmUnif} and \ref{thmUnifBV}.

We would like to mention
that Andr\'{e} Pierro de Camargo  suggested
another proof of the second conjecture for
functions with continuous third derivatives.
For odd $n$,
Theorem \ref{thm2} in Section \ref{secNum} indicates that
\[
\wfc{N_n}{f,x} - \wfc{f}{x} \wfc{D_n}{x} \approx
\frac{\wfc{f}{x} - \wfc{f}{1}}{2 \wlr{x - 1}} -
\frac{\wfc{f}{x} - \wfc{f}{-1}}{2 \wlr{x + 1}},
\]
and by solving this expression for $\wfc{f}{x}$ we derive the interpolant
\[
\wfc{f}{x} \approx \wfc{\tilde{B}_n}{f,x} :=
\frac{
\wfc{N_n}{f,x} + \frac{\wfc{f}{-1}}{2 \wlr{x + 1}}  - \frac{\wfc{f}{1}}{2 \wlr{x - 1}}
}
{
\wfc{D_n}{x} + \frac{1}{2 \wlr{x + 1}}  - \frac{1}{2 \wlr{x - 1}}
}.
\]
Note that $\tilde{B}_n$ is obtained by changing the absolute value of the
first and last weights of the interpolant in Equation \pRef{bary}
from $1$ to $1/2$. A similar argument applies to even $n$ and
the resulting barycentric interpolant $\tilde{B}_n$
has better convergence properties than Berrut's interpolant.
In fact, $\tilde{B}_n$ is
the interpolant corresponding to $d = 1$ in the
Floater-Hormann family \cite{Floater}, and
using the theory presented in \cite{Floater} we could
prove the second conjecture for
$f \in \wrm{C}^3$ by analyzing the asymptotic behavior
of $B_n - \tilde{B}_n$.

In summary, the present article shows that actually,
from the perspective of order of approximation,
Berrut's interpolants are biased by the functions
$\wfc{O}{f}$ and  $\wfc{E}{f}$, and we
see little reason for using them instead of
the interpolant $\tilde{B}_n$ above. In fact, in his
latter work \cite{Berrut2} prof. Berrut himself has mentioned that
using half integer weights at the endpoints
instead of $\pm 1$ leads to a better convergence rate.

Theorem \ref{thmUnif} shows
that the interpolant $\tilde{B}_n$ has
order of approximation $\wfc{o}{1/n}$,
and the most relevant questions in this subject are not the ones
raised by professors Mastroianni and Szabados, and which we
discuss in detail here. It is our opinion
that it is more important to understand how we should
choose the weights in the barycentric
interpolants in order to improve them, so that
we can justify the expensive $2n + 3$ divisions
per evaluation required by these interpolants.
This will be the subject of our next article
about barycentric interpolation.

\section{Proof of the first conjecture}
\label{secFirst}

In this section we prove Conjecture \ref{conj} by presenting
 $f \in \wrm{Lip \ 1}$ such that,
for
\pbDef{deftn}
t_n := 1/n
\hspace{1cm} \wrm{and} \hspace{1cm}
n_j := 2^{2^j},
\peDef{deftn}
we have
\pbTClaim{main}
\wfc{B_{n_j}}{f,t_{n_j}} - \wfc{f}{t_{n_j}} =  \wfc{B_{n_j}}{f,t_{n_j}}
\geq \frac{\wfc{\ln}{n_j}}{20 \, n_j}.
\peTClaim{main}
The function $f$ is given by
\pbDef{defg}
\wfc{f}{x} := \sum_{i = 100}^\infty \wfc{f_{n_i}}{x},
\peDef{defg}
for functions $f_m$ defined for $m$ such that $\sqrt{m}$ is an
integer multiple of $4$,
as follows:
\begin{eqnarray}
\pDefLabel{defgmA}
\wfc{f_m}{x} := 0 & \wrm{for} \ x <  \frac{1}{m} \ \ \wrm{or} \ \ x \geq \frac{\sqrt{m} - 3}{m},  \\
\pDefLabel{defgmB}
\wfc{f_m}{x} := x - \frac{1}{m}  & \wrm{for} \ \ \frac{1}{m} \leq  x < \frac{2}{m}, \\
\pDefLabel{defgmD}
\wfc{f_m}{x} := \frac{4 p + 3}{m} - x & \wrm{for} \ \ 0 \leq p \leq \frac{\sqrt{m}- 8}{4}
                 \ \ \wrm{and} \ \  \frac{4 p + 2}{m} \leq x < \frac{4 p + 4}{m}, \\
\pDefLabel{defgmC}
\wfc{f_m}{x} := x - \frac{4 p + 1}{m} & \wrm{for} \ \ 1 \leq p \leq \frac{\sqrt{m}- 8}{4}
                 \ \ \wrm{and} \ \  \frac{4 p}{m} \leq x < \frac{4 p + 2}{m}, \\
\pDefLabel{defgmE}
\wfc{f_m}{x} := x - \frac{\sqrt{m} - 3}{m} & \wrm{for} \ \
                \frac{\sqrt{m} - 4}{m} \leq x < \frac{\sqrt{m} - 3}{m}.
\end{eqnarray}

\pbFigureBH{figGm}{The functions fm}
\begin{picture}(400,130)(-15,10)

\put(-10, 80){\vector(1,0){340}}
\put( 0, 20){\vector(0,1){120}}
\put(-4,  120){\line(1,0){8}}
\put(-15,  118){$\frac{1}{m}$}
\put(-4,   40){\line(1,0){8}}
\put(-17,  38){$\frac{-1}{m}$}

\put(10,76){\line(0,1){8}}
\put( 5,65){$\frac{1}{m}$}
\put(10,105){$R_0$}

\put(10,80){\line(1,2){20}}

\put(30,76){\line(0,1){8}}
\put(25,65){$\frac{2}{m}$}
\put(30,120){\line(1,-2){40}}
\put(38,110){$F_0$}

\put(50,76){\line(0,2){8}}
\put(43,65){$\frac{3}{m}$}

\put(70,76){\line(0,1){8}}
\put(64,65){$\frac{4}{m}$}
\put(70,40){\line(1,2){40}}

\put(90,76){\line(0,1){8}}
\put(87,65){$\frac{5}{m}$}

\put(110,76){\line(0,1){8}}
\put(105,65){$\frac{6}{m}$}

\put(110,120){\line(1,-2){40}}
\put(73,128){$H_p = \wrm{hat}_p = R_p \bigcup F_p$}
\put(119,105){$F_p = \wrm{fall}_p$}
\put(50,105){$R_p = \wrm{raise}_p$}

\put(130,76){\line(0,1){8}}
\put(124,65){$\frac{7}{m}$}

\put(150,76){\line(0,1){8}}
\put(145,65){$\frac{8}{m}$}

\put(150,40){\line(1,2){40}}

\put(170,76){\line(0,1){8}}
\put(190,76){\line(0,1){8}}
\put(247,62){$\frac{\sqrt{m} - 5}{m}$}
\put(250,120){\line(-1,-2){40}}

\put(210,76){\line(0,1){8}}
\put(230,76){\line(0,1){8}}

\put(277,90){$\frac{\sqrt{m} - 4}{m}$}
\put(250,125){$H_{\frac{\sqrt{m} - 8}{4}}$}

\put(250,76){\line(0,1){8}}
\put(270,76){\line(0,1){8}}
\put(290,76){\line(0,1){8}}
\put(310,76){\line(0,1){8}}

\put(310,62){$\frac{\sqrt{m} - 3}{m}$}
\put(290,40){\line(1,2){20}}
\put(290,40){\line(-1,2){40}}
\put(303,47){$R_{\frac{\sqrt{m} - 4}{4}}$}

\end{picture}
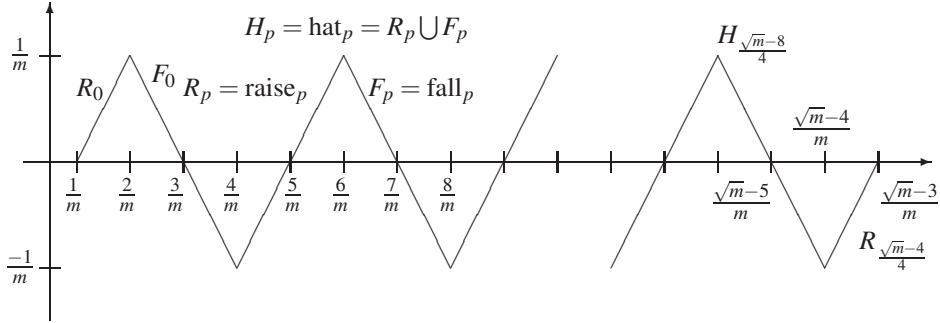
\caption{The function $f_m$. The support of $f_m$
is $[1/m,\wlr{\sqrt{m} - 3}/m]$.
The plot is divided in raise and fall regions, with
$R_p$ starting at $x = 4 p/m$ and $F_p$ starting at
$x = \wlr{4 p + 2}/m$. By joining $R_p$ and $F_p$ we
obtain the hat $H_p$.}
\peFigure{figGm}

Note that the series in Equation \pRef{defg} converges to
$f \in \wrm{Lip \ 1}$ because $n_j = 2^{2^j}$ and the identities
\pbTClaim{support}
n_j^2 = 2^{2^{j+1}} = n_{j+1} \Rightarrow
\frac{\sqrt{n_{j+1}} - 3}{n_{j+1}} <
\frac{1}{\sqrt{n_{j+1}}} = \frac{1}{n_j}
\peTClaim{support}
imply that the support of the functions $f_{n_j}$
are disjoint, and $f_m \in \wrm{Lip \ 1}$.

Equation \pRef{main} follows from Equation \pRef{defg} and the following Lemmas:

\pbLemmaBT{lemHeader}{The error for the first terms}
If $100 \leq i < j$ then
\pbTClaim{cheader}
\wfc{f_{n_i}}{t_{n_j}} = 0
\hspace{1cm} \wrm{and} \hspace{1cm}
\wfc{B_{n_j}}{f_{n_i},t_{n_j}} \geq -\frac{9}{8 n_j}.
\peTClaim{cheader}
\peLemma{lemHeader}

\pbLemmaBT{lemMain}{The error for the main term}
For $j \geq  100$ we have that
\pbTClaim{cmain}
\wfc{f_{n_j}}{t_{n_j}} = 0 \hspace{1cm} \wrm{and} \hspace{1cm}
\wfc{B_{n_j}}{f_{n_j},t_{n_j}} \geq \frac{\wfc{\ln}{n_j}}{16 n_j}.
\peTClaim{cmain}
\peLemma{lemMain}

\pbLemmaBT{lemTail}{The error for the last terms}
For $i > j \geq 100$ and $0 \leq k \leq n_j$ we have
\pbTClaim{ctail}
\wfc{f_{n_i}}{t_{n_j}} = 0
\hspace{1cm} \wrm{and} \hspace{1cm}
\wfc{f_{n_i}}{x_{k,n_j}} =  0.
\peTClaim{cTail}
\peLemma{lemTail}

The lemmas above show that $\wfc{f}{t_{n_j}} = 0$ for $j \geq 100$,
and the second part of Equation \pRef{main} follows from these lemmas because
\[
\wfc{B_{n_j}}{f,t_{n_j}} =
\left.
\sum_{k = 0}^{n_j}
 \wlr{-1}^k  \frac{\sum_{i = 100}^\infty \wfc{f_{n_i}}{x_{k,n_j}} }{t_{n_j} - x_{k,n_j}}
\right/
\sum_{k=0}^{n_j} \wlr{-1}^k  \frac{1}{t_{n_j} - x_{k,n_j}}
\]
\[
=
\left.
\sum_{k = 0}^{n_j}
\wlr{-1}^k  \frac{\sum_{i = 100}^{j} \wfc{f_{n_i}}{x_{k,n_j}}}{t_{n_j} - x_{k,n_j}}
\right/
\sum_{k = 0}^{n_j} \wlr{-1}^k  \frac{1}{t_{n_j} - x_{k,n_j}}
\]
\[
=
\sum_{i = 100}^{j}
\left.
\sum_{k = 0}^{n_j}
\wlr{-1}^k  \frac{\wfc{f_{n_i}}{x_{k,n_j}}}{t_{n_j} - x_{k,n_j}}
\right/
 \sum_{i=0}^n \wlr{-1}^k  \frac{1}{t_{n_j} - x_{k,n_j}}
\]
\pbTClaim{final}
= \sum_{i = 100}^{j} \wfc{B_{n_j}}{f_{n_i}, t_{n_j}}
\geq \wfc{B_{n_j}}{f_{n_j}, t_{n_j}} - 9 \frac{j - 100} {8 n_j}
 \geq \frac{ \wfc{\ln}{n_j}}{16 \, n_j} - 9 \frac{j} {8 n_j}
\peTClaim{final}
and, finally, the reader can verify that for $j \geq 100$
\[
j < \wfc{\ln}{2^{2^j}}/1000 = \wfc{\ln}{n_j}/1000
\]
and Equation \pRef{main} follows from Equation \pRef{final}.

We end this section presenting a proof of
the lemmas above and one more lemma:



\pbLemmaBT{lemHarmo}{Shifted harmonic sums}
If $a > 0$ and $\ell \geq 1$ is an integer then
\pbTClaim{harmo}
\sum_{j = 0}^{\ell-1} \frac{1}{a + j} \geq \wfc{\ln}{a + \ell} - \wfc{\ln}{a} + \frac{1}{2 a} - \frac{1}{2 \wlr{a + \ell}}.
\peTClaim{harmo}
\peLemma{lemHarmo}


\pbProofB{Lemma}{lemHeader}
If $i < j$ then $n_i < n_j$,
$t_{n_j} = 1 / n_j < 1/n_i$
and Equation \pRef{defgmA} implies that $\wfc{f_{n_i}}{t_{n_j}} = 0$.
This proves the first part of Equation \pRef{cheader}.
Let $\wfc{N_{n_j}}{f_{n_i},t_{n_j}}$ and
$\wfc{D_{n_j}}{t_{n_j}}$ be as in Equation \pRef{numDen}.
Lemma \ref{lemDen} in Section \ref{secDen}
shows that $\wabs{\wfc{D_{n_j}}{t_{n_j}}} \geq n_j$
and this reduces the proof of Lemma \ref{lemHeader} to
the verification of the equation
\pbTClaim{claimND}
\wfc{N_{n_j}}{f_{n_i},t_{n_j}} \geq - 3/4,
\peTClaim{claimND}
as we do below.
Note that the definition of $n_j$ in Equation \pRef{deftn} implies that
if $100 \leq i < j$ then, for $m = n_i$,
\[
n_j = 4 q m
\hspace{0.5cm}
\wrm{with}
\hspace{0.5cm}
q \geq 16,
\hspace{0.5cm}
t_{n_j} = \frac{1}{n_j} = \frac{1}{4 q m}
\hspace{1cm} \wrm{and} \hspace{1cm}
x_{k,n_j} = \frac{2k}{4 q m } - 1.
\]
Equation \pRef{defgmA} shows that
\[
\wfc{f_m}{x_{k,n_j}} = \wfc{f_m}{\frac{2 k}{4 q m} - 1} = 0
\hspace{0.5cm} \wrm{if} \hspace{0.5cm}
k < 2 q m + 2 q \hspace{0.5cm} \wrm{or} \hspace{0.5cm}
k \geq 2 q m + 2 q \wlr{\sqrt{m} - 3},
\]
and Equation \pRef{numDen}, with the index $k$ replaced by $k + 2 q m$,
leads to
\[
\wfc{N_{n_j}}{f_{n_i},t_{n_j}} = \sum_{k = 2 q}^{2 q \wlr{\sqrt{m} - 3} - 1}
 \wlr{-1}^{k} \frac{\wfc{f_{m}}{\frac{2k}{4 q m}} }{\frac{1}{4 q m} - \frac{2k}{4q m }}
= 4 q m \sum_{k = 2 q}^{2 q \wlr{\sqrt{m} - 3} - 1} \wlr{-1}^{k + 1}
\frac{\wfc{f_{m}}{\frac{2k}{4 q m}}  }{2 k - 1}.
\]
Motivated by Figure \ref{figGm},
we split the parcels of $\wfc{N_{n_j}}{f_{n_i},t_{n_j}}$ in
$h := \wlr{\sqrt{m} - 8}/4$ hats plus the last half of
$R_0$, which we call by $R_-$,
the part $F_0$, and the first half of $R_{\wlr{\sqrt{m} - 4}/4}$,
which we call by $R_+$. Formally we have
\[
\frac{\wfc{N_{n_j}}{f_{n_i},t_{n_j}}}{4 q m} = R_- + F_0 +
\wlr{\, \sum_{p = 1}^{h} \wlr{F_p + R_p} \,} + R_+
= R_- + F_0 + \wlr{\, \sum_{p = 1}^{h} H_p } + R_+,
\]
for
\begin{eqnarray}
\pDefLabel{uminus}
R_- & := & \sum_{k = 2 q}^{4 q - 1}  \wlr{-1}^{k + 1}
\frac{\wfc{f_{m}}{\frac{2k}{4 q m}} }{2 k - 1}, \\
\pDefLabel{fallp}
F_p & := & \sum_{k = 8 p q + 4 q}^{8 \wlr{p + 1} q - 1}
 \wlr{-1}^{k + 1}  \frac{\wfc{f_{m}}{\frac{2k}{4 q m}}}{2 k - 1}, \\
\pDefLabel{uj}
R_p & := & \sum_{k = 8 p q}^{8 p q + 4 q - 1} \wlr{-1}^{k + 1}
\frac{\wfc{f_{m}}{\frac{2k}{4 q m}} }{2 k - 1}, \\
\pDefLabel{uplus}
R_+ & := & \sum_{k = 2 q \wlr{2 \sqrt{m} - 4}}^{q \wlr{\sqrt{m} - 3} - 1} \wlr{-1}^{k + 1}
\frac{\wfc{f_{m}}{\frac{2k}{4 q m}}  }{2 k - 1}, \\
\nonumber
H_p & :=  &  R_p + F_p,
\end{eqnarray}
and to prove Equation \pRef{claimND}
it suffices to show that $R_-, R_+, H_p > 0$ and
\pbTClaim{boundD0}
F_0 \geq \frac{-3}{16 q m},
\peTClaim{boundD0}
and this is done from this point to the end of this
proof.

Let us start by writing $R_p$ as a sum of positive terms.
In raising ranges $f_m$ is defined by Equations \pRef{defgmB},
\pRef{defgmC} and \pRef{defgmE},
and Equation \pRef{uj} yields
\[
R_p = \sum_{k = 8 p q}^{8 p q + 4 q - 1} \wlr{-1}^{k + 1}
\frac{\wlr{\frac{2k}{4 q m} - \frac{4 p + 1}{m}}  }{2 k - 1}
\]
\[
= \frac{1}{4 q m} \sum_{k = 8 p q}^{8 p q + 4 q - 1} \wlr{-1}^{k + 1}
  \frac{\wlr{2 k - 4 q - 16 p q }  }{2 k - 1}.
\]
Splitting the indexes $k$ in even and odd groups we obtain
\[
R_p = \frac{-1}{4 q m} \sum_{\ell = 4 p q}^{4 p q + 2 q - 1}
\wlr{
\frac{4 \ell - 4 q - 16 p q}{4 \ell - 1} -
\frac{4 \ell + 2 - 4 q - 16 p q}{4 \ell + 1}
}
\]
\[
=
\frac{-1}{4 q m} \sum_{\ell = 4 p q}^{4 p q + 2 q - 1}
\wlr{
\frac{1 - 4 q - 16 p q}{4 \ell - 1} -
\frac{1 - 4 q - 16 p q}{4 \ell + 1}
},
\]
and
\pbTClaim{boundUj}
R_p =
\frac{16 p q + 4 q - 1}{2 q m}
\sum_{\ell = 4 p q}^{4 p q + 2 q - 1} \frac{1}{16 \ell^2 - 1} > 0.
\peTClaim{boundUj}
The same argument
using Equations \pRef{uminus} and \pRef{uplus}
shows that $R_-,R_+ > 0$.
Similarly, for $F_p$ Equations \pRef{defgmD} and \pRef{fallp} lead to
\[
F_p = \sum_{k =  8 p q + 4 q}^{8 \wlr{p+1} q - 1}
\wlr{-1}^{k + 1} \frac{\wlr{\frac{4 p + 3}{m} - \frac{2k}{4 q m}} }{2 k - 1}
\]
\[
= \frac{1}{4 q m} \sum_{k =  8 p q + 4 q}^{8 \wlr{j+1} q - 1}
\wlr{-1}^{k + 1}  \frac{\wlr{16 p q + 12 q - 2 k} }{2 k - 1}.
\]
As before,
\[
F_p = \frac{1}{4 q m} \sum_{\ell = 4 p q + 2 q}^{4 \wlr{p+1} q - 1}
\wlr{
\frac{4 \ell - 16 p q - 12 q}{4 \ell - 1} -
\frac{4 \ell + 2 - 16 p q - 12 q}{4 \ell + 1}
}
\]
\[
=  \frac{1}{4 q m} \sum_{\ell = 4 p q + 2 q}^{4 \wlr{p+1} q - 1}
\wlr{
\frac{1 - 16 p q - 12 q}{4 \ell - 1} -
\frac{1 - 16 p q - 12 q}{4 \ell + 1}
},
\]
and
\pbTClaim{boundDj}
F_p = - \frac{16 p q + 12 q - 1}{2 q m} \sum_{\ell = 4 p q + 2 q}^{4 \wlr{p+1} q - 1}
\frac{1}{16 \ell^2 - 1}.
\peTClaim{boundDj}
In particular, for $p = 0$ we have
\[
F_0 = - \frac{12 q - 1}{2 q m} \sum_{\ell = 2 q}^{4 q - 1}
\frac{1}{16 \ell^2 - 1}
\geq - \frac{12 q - 1}{2 q m} \frac{2 q}{64 q^2 - 1}
= - \frac{1}{2 q m} \frac{24 q^2 - 2 q}{64 q^2 - 1} \geq \frac{-3}{16 q m},
\]
and this proves Equation \pRef{boundD0}.

We now show that, for $p \geq 1$, $H_p = R_p + F_p > 0$.
Replacing $\ell$ by $k + 2 q$ in Equation \pRef{boundDj}
and $\ell$ by $k$ in Equation \pRef{boundUj}
we obtain
\[
H_p = \frac{1}{2q m} \sum_{k = 4 p q}^{4 p q + 2 q - 1} a_k,
\]
for
\[
a_k = \frac{16 p q + 4 q - 1}{16 k^2 - 1}
    - \frac{16 p q + 12 q - 1}{16 \wlr{k + 2 q}^2 - 1},
\]
and our final goal is to show that $a_k > 0$.
We can write $a_k$ as $u_k/v_k$ for
\[
u_k := \wlr{16 p q +  4 q - 1} \wlr{16 \wlr{k + 2 q}^2 - 1} -
       \wlr{16 p q + 12 q - 1} \wlr{16 k^2 - 1}
\]
and
\[
v_k := \wlr{16 k^2 - 1} \wlr{16 \wlr{k + 2 q}^2 - 1}.
\]
The denominator $v_k$ is clearly positive, and in order
to analyze $u_k$ we replaced $k$ by $4 p q + \xi q$,
with $\xi \in [0,2)$,
and used Wolfram Alpha to obtain
\[
u_k = 8 q \wlr{256 p^2 q^2 + 256 p q^2 - 32 p q - 16 q^2 \xi^2 + 32 q^2 + 32 q^2 \xi - 8 q \xi - 8 q + 1}.
\]
Since we are concerned with $q \geq 16$, $p \geq 1$ and $\xi \in [0,2)$, it is
clear that $u_k > 0$ and
the proof of Lemma \ref{lemHeader} is complete.
\peProof{Lemma}{lemHeader}\\


\pbProofB{Lemma}{lemMain} Let us write $m = n_j$.
According to Equation \pRef{deftn}, $t_m = 1 / m$,
and Equation \pRef{defgmB} yields $\wfc{f_{m}}{t_{m}} = 0$.
We have that
\[
\wfc{B_m}{f_m,t_m} = \wfc{N_m}{f_m,t_m} / \wfc{D_m}{t_m}
\]
for $\wfc{N_m}{f_m,t_m}$
and $\wfc{D_m}{t_m}$ in Equation \pRef{numDen}.
Since
\[
 \frac{1}{2m} = x_{\frac{m}{2},m} + \frac{1}{2m}
 < t_m = \frac{1}{m} < x_{\frac{m}{2} + 1,m} - \frac{1}{2m} = \frac{3}{2m}
\]
and $m$ is a multiple of four and we have that
\[
t_m = \frac{2 \times \wlr{2 m} + 0 + 1}{4m} - 1.
\]
Equation \pRef{boundDen} in Section \ref{secDen} with
$\wfc{\rho_n}{x} = 0$ shows that
\[
0 < \wfc{D_m}{t_m} \leq \wfc{A}{0} + 1/2 = \pi m / 2 + 1/2 < 4 m,
\]
and in order to prove Lemma \ref{lemMain} it
suffices to show that
\pbTClaim{mainR}
\wfc{N_m}{f_m,t_m} \geq \wfc{\ln}{m}/ 4.
\peTClaim{mainR}
This is our goal now.
Equations \pRef{defgmA}--\pRef{defgmE} imply
that $\wfc{f_m}{2 k/m - 1} = \wlr{-1}^{k+1} / m$ for
\[
k = m/2 + 1, \dots,  m/2 + \wlr{\sqrt{m} - 4}/2
\]
and $\wfc{f_m}{2 k / m - 1} = 0$ for the remaining $k$s
(see Figure \ref{figGm}.) Therefore,
\[
\wfc{N_m}{f_m,t_m} := \sum_{k = m/2 + 1}^{m/2 + \frac{\sqrt{m} - 4}{2}}
\wlr{-1}^{k}  \frac{\wfc{f_{m}}{2 k/m - 1}}{1/m - 2k/m + 1}.
\]
Making the change of indexes $k = m/2 + i$ and noting that $m/2$ is even
we obtain
\[
\wfc{N_m}{f_m,t_m} = \sum_{i = 1}^{\frac{\sqrt{m} - 4}{2}}
 \wlr{-1}^{i}  \frac{\wlr{-1}^{i + 1}/m}{1/m - 2 i/m}
= \sum_{i = 1}^{\frac{\sqrt{m} - 4}{2}} \frac{1}{2 i - 1} =
  \frac{1}{2} \sum_{i = 0}^{\frac{\sqrt{m} - 6}{2}} \frac{1}{i + 1/2},
\]
and Lemma \ref{lemHarmo} with $a = 1/2$ and $\ell = \wlr{\sqrt{m} - 4}/2$ yields
\[
2 \wfc{N_m}{f_m,t_m} \geq
\wfc{\ln}{\frac{\sqrt{m} - 3}{2}} - \wfc{\ln}{1/2} + 1 -
\frac{1}{2 \wlr{\frac{1}{2} + \frac{\sqrt{m} - 4}{2}}}
=
\wfc{\ln}{\sqrt{m} - 3} + 1 - \frac{1}{\sqrt{m} - 3}.
\]
Therefore,
\pbTClaim{mainR2}
\wfc{N_m}{f_m,t_m} = \wfc{\ln}{\sqrt{m}}/2 +  \delta_m / 2 =  \wfc{\ln}{m}/4 + \delta_m / 2
\peTClaim{mainR2}
for
\[
\delta_m = 1 + \wfc{\ln}{1 - \frac{3}{\sqrt{m}}} - \frac{1}{\sqrt{m} - 3}.
\]
Since $4 m \geq 2^{2^{100}}$ we have that $\delta_m > 0$. Therefore,
Equation \pRef{mainR2} implies Equation \pRef{mainR} and this
proof is complete.
\peProof{Lemma}{lemMain}\\


\pbProofB{Lemma}{lemTail}
Equation \pRef{support} implies that if $i > j$ then
\[
\frac{\sqrt{n_i} - 3}{n_i} < \frac{1}{\sqrt{n_i}} \leq \frac{1}{n_j} = t_{n_j},
\]
and Equation \pRef{defgmA} implies that $\wfc{f_{n_i}}{t_{n_j}} = 0$.
In order to show that $\wfc{f_{n_i}}{x_{k,n_j}} = 0$ we
recall that $x_{k,n_j} = 2 k / n_j - 1$ and analyze
two possibilities:
\begin{itemize}
\item[(i)]
If $k \leq n_j/2$ then $x_{k,n_j} \leq 0$, and Equation \pRef{defgmA}
implies that $\wfc{f_{n_i}}{x_{k,n_j}} = 0$.
\item[(ii)] If $k > n_j/2$ then
\[
x_{k,n_j} \geq \frac{2}{n_j} \geq \frac{2}{\sqrt{n_i}} > \frac{\sqrt{n_i} - 3}{n_i},
\]
and Equation \pRef{defgmA} shows that $\wfc{f_{n_i}}{x_{k,n_j}} = 0$.
\end{itemize}
Therefore, $\wfc{f_{n_i}}{x_{k,n_j}} = 0$ in both cases
we have proved Lemma \ref{lemTail}.
\peProof{Lemma}{lemTail}\\


\pbProofB{Lemma}{lemHarmo}
For $b > 0$, let $h_b: [0,1] \rightarrow \wrone{}$ be the function
\[
\wfc{h_b}{t} := \frac{1}{b} + \wlr{\frac{1}{b + 1} - \frac{1}{b}} t - \frac{1}{b + t}
= \frac{1}{b} - \frac{t}{b \wlr{b + 1}} - \frac{1}{b + t}.
\]
Since $\wfc{h_b}{0} = \wfc{b_b}{1} = 0$ and $h_b$ is concave
we have that $h_b \geq 0$.
Therefore,
\[
0 \leq \int_0^1 \wfc{h_b}{t} dt = \frac{1}{b} - \frac{1}{2 b \wlr{b + 1}} - \wfc{\ln}{b + 1} + \wfc{\ln}{b}
\]
and
\[
\frac{1}{b} \geq  \frac{1}{2 b \wlr{b + 1}} + \wfc{\ln}{b + 1} - \wfc{\ln}{b} =
\frac{1}{2b} - \frac{1}{ 2\wlr{b + 1}} + \wfc{\ln}{b + 1} + \wfc{\ln}{b}.
\]
It follows that
\[
\sum_{j = 0}^{\ell - 1} \frac{1}{a + j} \geq
\sum_{j = 0}^{\ell - 1}
\frac{1}{2 \wlr{a + j}} - \frac{1}{2 {\wlr{a + j + 1}}} + \wfc{\ln}{a + j + 1} - \wfc{\ln}{a + j} =
\]
\[
\wfc{\ln}{a +\ell} - \wfc{\ln}{a} + \frac{1}{2 a} - \frac{1}{2 \wlr{a + \ell}}
\]
and we are done.
\peProof{Lemma}{lemHarmo}\\

\section{The denominator}
\label{secDen}

In this section we analyse the denominator $\wfc{D_n}{x}$ of the
interpolant in Equation \pRef{bary}, using the
function $A: [0,1) \rightarrow \wrone{}$ given by
\pbDef{defA}
\wfc{A}{x} := \sum_{k = 0}^\infty \wlr{-1}^k \frac{4 k + 2}{\wlr{2 k + 1}^2 - x}.
\peDef{defA}
This function is increasing and
can be extended to a homeomorphism between $[0,1]$ and
$[\pi/2,+\infty] \subset \overline{\wrone{}}$,
with the topology in the introduction,
as shown by the next lemma. In the
rest of the article we work with
this extension of $A$ and its inverse $A^{-1}$.

\pbLemmaBT{lemA}{The function $A$}
The function $A$ defined in Equation \pRef{defA} is increasing,
$\wfc{A}{0} = \pi/2$ and
\pbTClaim{boundA}
 -1/2 \leq \wfc{A}{x} - \frac{2}{1 - x} \leq \frac{\pi - 4}{2} < -0.42.
\peTClaim{boundA}
In particular, $A$ can be extended to a homeomorphism
between $[0,1]$ and $[\pi/2,+\infty]$.
\peLemma{lemA}

The section is based upon the observation that
for a regular $x$, as $j$ tends to infinity the
denominator $\wfc{D_{n_j}}{x}$  can be accurately described by
the expression
\[
\wfc{D_{n_j}}{x} \approx \wlr{-1}^{\wfc{\iota_{n_j}}{x}} \, {n_j} \ \wfc{A}{\wfc{\rho_{n_j}^2}{x}},
\]
where $A$ is the function defined in Equation \pRef{defA},
\pbDef{rho}
\wfc{\iota_n}{x} := \wfloor{n \wlr{x + 1}/2}
\hspace{1cm} \wrm{and} \hspace{1cm}
\wfc{\rho_n}{x} := n \wlr{x - x_{\wfc{\iota_n}{x},n}} - 1,
\peDef{rho}
so that $\wfc{\iota_n}{x} \in \wset{0,\dots, n - 1}$,

\pbDef{rhoDec}
x = \frac{2 \wfc{\iota_n}{x} + \wfc{\rho_n}{x} + 1}{n} - 1
\hspace{1cm} \wrm{and} \hspace{1cm}
\wfc{\rho_n}{x} \in \wlr{-1,1}.
\peDef{rhoDec}

Formally, we have the following lemma:
\pbLemmaBT{lemDen}{The size and sign of the denominator}
If $x \in (-1,1) \setminus \wset{x_{0,n},\dots,x_{n,n}}$ then
\[
\wfc{\rho_n}{x} \in (-1,1), 
\hspace{1cm} \wfc{\wrm{sign}}{\wfc{D_n}{x}} = \wlr{-1}^{\wfc{\iota_n}{x}},
\]
\pbTClaim{boundDen}
 \wabs{ \wabs{ \wfc{D_n}{x} / n} - \wfc{A}{\wfc{\rho_n^2}{x}}} \leq
 \frac{1}{4 \wlr{1 + \wfc{\iota_n}{x}}} + \frac{1}{4 \wlr{n - \wfc{\iota_n}{x}}}
 \leq \frac{1}{2},
\peTClaim{boundDen}
and
\pbTClaim{boundDenA}
 \wabs{ \wfc{D_n}{x} / n} \geq 1
 \hspace{0.5cm} \wrm{and} \hspace{1cm}
 \wabs{ \wfc{D_n}{x} / n} \geq \wfc{A}{\wfc{\rho_n^2}{x}}/2  \geq
 \frac{3}{4 \wlr{1 - \wfc{\rho_n^2}{x}}}.
\peTClaim{boundDenA}
In particular, if $x$ is regular for the sequence $n_j$ then
\pbTClaim{boundDenR}
\lim_{j \rightarrow \infty} \frac{\wabs{\wfc{D_{n_j}}{x}}}{n_j \, \wfc{A}{\wfc{\rho_{n_j}^2}{x}}} = 1.
\peTClaim{boundDenR}
\peLemma{lemDen}

The last two lemmas imply that the possible values for
$\lim_{j \rightarrow \infty} \wfc{D_{n_j}}{x}/n_j$
can be found by analysing the limits
$\lim_{j \rightarrow \infty} \wfc{\iota_{n_j}}{x}$ and
$\lim_{j \rightarrow \infty} \wfc{\rho_{n_j}^2}{x}$.

\pbCorolBT{corConv}{Convergence of $\wfc{D_n}{x}/n$}
If $x$ is regular for the sequence $n_j$ then
\pbTClaim{corConvLhs}
\lim_{j \rightarrow \infty} \frac{1}{n_j} \wfc{D_{n_j}}{x} = L \in \overline{\wrone{}}
\peTClaim{corConvLhs}
if and only if
\pbTClaim{corConvRhs}
\lim_{j \rightarrow \infty} \wlr{-1}^{\wfc{\iota_{n_j}}{x}} = \wfc{\wrm{sign}}{L},
\hspace{0.7cm}
\wabs{L} \geq \frac{\pi}{2}
\hspace{0.7cm} \wrm{and} \hspace{0.7cm}
\lim_{j \rightarrow \infty} \wfc{\rho_{n_j}^2}{x} = \wfc{A^{-1}}{\wabs{L}}.
\peTClaim{corConvRhs}
\peCorol{corConv}

This corollary leads to a clean description of the
limits $\lim_{j \rightarrow \infty} \wfc{D_{n_j}}{x}/n_j$
when $x$ is irrational, due to the following theorem by
S. Hartmann:

\pbTheoremBT{thmH}{Hartmann's Theorem \cite{Hartmann}}
For every irrational number $\xi$, and integers $s,a,b$, with $s \geq 1$,
there are infinitely many integers $u$ and $v > 0$ such that
\pbDef{hart}
\wabs{\xi - \frac{u}{v}} \leq \frac{2 s^2}{v^2}
\hspace{0.5cm} \wrm{with} \hspace{0.5cm}
 u \equiv a \ \wrm{mod} \ s
 \hspace{0.5cm} \wrm{and} \hspace{0.5cm}
 v \equiv b \ \wrm{mod} \  s.
\peDef{hart}
\peTheorem{thmH}

Using Hartmann's theorem we can prove the following Lemma:

\pbLemmaBT{lemIrrat}{Convergence of the denominator for irrational $x$}
If $x \in (-1,1)$ is irrational then for each $r_n \in \wset{0,1}$ and
$y$ with $\wabs{y} \in [\pi/2,+\infty]$
there exists a sequence $n_j$ such that
\pbTClaim{limIrrat}
n_j \equiv r_n \ \wrm{mod} \ 2
\hspace{0.7cm} \wrm{and} \hspace{0.7cm}
\lim_{j \rightarrow \infty} \frac{1}{n_j} \wfc{D_{n_j}}{x} = y.
\peTClaim{limIrrat}
\peLemma{lemIrrat}

In words, Lemma \ref{lemIrrat} shows that if $x$ is irrational
then we can obtain all elements in the extended intervals
$[-\infty, -\pi/2] $ and $[\pi/2,+\infty]$
as limits for $\wfc{D_{n_j}}{x}/n_j$, for sequences
$n_j$ with the same parity, be this parity odd or even. Unfortunately
things are more complex when $x$ is rational
and we must consider a few cases, as we do in
the next lemmas. The first one shows that the set of possible
limits for $\wfc{D_n}{x}/n$ is finite in this case.

\pbLemmaBT{lemRat}{Finitely many limits $\wfc{D_n}{x}/n$ for $x$ rational}
For $p,q \in \wn{}$, with $q \neq 0$.
If $x = p/q - 1 \in \wlr{-1,1}$  is regular for the sequence $n_j$
and
\pbDef{limRat}
\lim_{j \rightarrow \infty} \frac{1}{n_j} \wfc{D_{n_j}}{x}
= L \in \overline{\wrone{}}
\peDef{limRat}
then $L$ is finite and $\wabs{L} = \wfc{A}{m^2/q^2}$
for some $m \in \wz{}$  with $0 \leq m < q$. Moreover,
there exists $j_0$ such that if $j \geq j_0$ then
\[
\wlr{-1}^{\wfc{\iota_{n_j}}{x}} = \wfc{\wrm{sign}}{L}
\hspace{1cm} \wrm{and} \hspace{1cm}
\wabs{\wfc{\rho_{n_j}}{x}} = m / q.
\]
\peLemma{lemRat}

The hypothesis of Lemma \ref{lemRat}
accounts for $L = \pm \infty$, but
its thesis states that this case is actually not
possible. In particular, this Lemma implies that if
$x$ is regular for $n_j$ then the
$\wfc{D_{n_j}}{x}/n_j$ are bounded.
Lemma \ref{lemRat} also shows that
if the sequence $\wfc{D_{n_j}}{x}/{n_j}$ converges
in $\overline{\wrone{}}$ and $x$ is rational
and regular
then $\wabs{\wfc{\rho_{n_j}}{x}}$ and the parity of
$\wfc{\iota_{n_j}}{x}$ become eventually
constant, and $L$ belongs to one of the
two finite sets
\begin{eqnarray}
\nonumber{}
\wfc{\wcal{O}}{p/q} := \left\{
L \in \overline{\wrone{}} \ \ \wrm{such \ that \ Equation
\ \pRef{limRat} \ holds \ for \ some \ odd \ sequence \ } n_j \
\right. \\
\pDefLabel{oddpq}
\left.
\wrm{for \ which } \ x = p/q - 1 \ \wrm{is \ regular }  \right\}
\end{eqnarray}
and
\begin{eqnarray}
\nonumber{}
\wfc{\wcal{E}}{p/q} := \left\{
L \in \overline{\wrone{}} \ \ \wrm{such \ that \ Equation
\ \pRef{limRat} \ holds \ for \ some \ even \ sequence \ } n_j \
\right. \\
\pDefLabel{evenpq}
\left.
\wrm{for \ which } \ x = p/q - 1 \ \wrm{is \ regular }  \right\}.
\end{eqnarray}

The description of the sets of limits
$\wfc{\wcal{O}}{p/q}$ and $\wfc{\wcal{E}}{p/q}$
is a tedious exercise in elementary
number theory, but we present it below
for completeness. The possible cases
are listed in the next three corollaries.
After the statement of these corollaries we end
this section  with the proofs of the result stated in it.
\pbCorolBT{corRatPQ}{$\wfc{\wcal{O}}{p/q}$ and $\wfc{\wcal{E}}{p/q}$  for odd $p$ and $q$}
If $\wfc{\wrm{gcd}}{p,q} = 1$  and $p$ and $q$ are odd then
\pbTClaim{oddRatPQ}
\wfc{\wcal{O}}{p/q} =
\wset{\pm \wfc{A}{4 \ell^2/q^2} \ \wrm{with}  \ \ell \in \wset{0,1,\dots, \wlr{q - 1}/2}}
\peTClaim{oddRatPQ}
and
\pbTClaim{evenRatPQ}
\wfc{\wcal{E}}{p/q} =
\wset{\pm \wfc{A}{\wlr{2 \ell + 1}^2/q^2} \ \wrm{with} \ \ell \in \wset{0,1,\dots, \wlr{q - 3}/2}}.
\peTClaim{evenRatPQ}
\peCorol{corRatPQ}

\pbCorolBT{corRat2PQ}{$\wfc{\wcal{O}}{2p/q}$ and $\wfc{\wcal{E}}{2p/q}$}
If $\wfc{\wrm{gcd}}{p,q} = 1$ and $q$ is odd then
\begin{eqnarray}
\nonumber{}
\wfc{\wcal{O}}{2p/q} =
\left\{ \wlr{-1}^s \wfc{A}{\wlr{4 \ell + 2 p - 2 s - q}^2/q^2} \
\wrm{for} \ s \in \wset{0,1} \ \ \wrm{and} \ \ \right.\\
\pDefLabel{oddRat2PQ}
\left.
\ell \in \wz{} \ \ \wrm{with} \ \
s - p + 1 \leq 2 \ell \leq s - p + q - 1
\right\}
\end{eqnarray}
and
\begin{eqnarray}
\nonumber{}
\wfc{\wcal{E}}{2p/q} =
\left\{ \wlr{-1}^s \wfc{A}{\wlr{4 \ell - 2 s - q}^2/q^2} \
\wrm{for} \ s \in \wset{0,1} \ \ \wrm{and} \ \ \right.\\
\pDefLabel{evenRat2PQ}
\left.
\ell \in \wz{} \ \ \wrm{with} \ \
s + 1 \leq 2 \ell \leq s + q - 1
\right\}.
\end{eqnarray}
\peCorol{corRat2PQ}

Finally,
\pbCorolBT{corRatP2Q}{$\wfc{\wcal{O}}{p/2q}$ and $\wfc{\wcal{E}}{p/2q}$}
If $\wfc{\wrm{gcd}}{p,q} = 1$ and $p$ is odd then
\pbTClaim{oddRatP2Q}
\wfc{\wcal{O}}{p/2q} =
\wset{ \pm \wfc{A}{\frac{\wlr{2 \ell + 1}^2}{4q^2}}  \ \wrm{with}  \ \ell \in \wset{0,1,\dots, q - 1}}
\peTClaim{oddRatP2Q}
and
\pbTClaim{evenRatP2Q}
\wfc{\wcal{E}}{p/2q} =
\wset{ \pm \wfc{A}{\ell^2/q^2}  \ \wrm{with}  \ \ell \in \wset{0,1,\dots, q - 1}}.
\peTClaim{evenRatP2Q}
\peCorol{corRatP2Q}


\pbProofB{Lemma}{lemA}
The derivative of $A$
\[
\wdfc{A}{x} = \sum_{k = 0}^\infty \wlr{-1}^k \frac{4 k + 2}{\wlr{\wlr{2 k + 1}^2 - x}^2}
\]
has parcels of alternating
signs and decreasing absolute values, with a positive first term.
Therefore $\wdfc{A}{x} > 0$ for all $x \in [0,1)$, and $A$ is a increasing
function of $x$.  Moreover, executing the command
\begin{verbatim}
Sum[ 2 (-1)^k / (2 k + 1), k = 0 to Infinity ]
\end{verbatim}
in the software Wolfram Alpha we obtain that $\wfc{A}{0} = \pi/2$.

The same argument used above shows that the function
$h: [0,1] \rightarrow \wrone{}$ given by
\[
\wfc{h}{x} = \wfc{A}{x} - \frac{2}{1 - x} = - \sum_{k = 1}^\infty
\wlr{-1}^k \frac{4 k + 2}{\wlr{2 k + 1}^2 - x}
\]
is decreasing, and executing
\begin{verbatim}
Sum[ 2 (-1)^k / (2 k + 1), k = 1 to Infinity ]
\end{verbatim}
and
\begin{verbatim}
Sum[ (-1)^k (4 k + 2)/ ((2 k + 1)^2 - 1), k = 1 to Infinity ]
\end{verbatim}
in Wolfram Alpha we obtain that
\[
\wfc{h}{0} =  \wlr{\pi - 4}/2 \approx -0.429036
\hspace{1cm} \wrm{and} \hspace{1cm}
\wfc{h}{1} =  -1/2.
\]
This proves Equation \pRef{boundA}.
\peProof{Lemma}{lemA}\\


\pbProofB{Lemma}{lemDen}
We have that
\[
x_{\wfc{\iota_n}{x},n} < x < x_{\wfc{\iota_n}{x} + 1,n}
\ \
\Rightarrow
\ \
0 < \wfc{\theta_n}{x} := n \wlr{x + 1}/2 - \wfc{\iota_n}{x} < 1.
\]
Equation \pRef{rho} defines $\wfc{\rho_n}{x} := n \wlr{x -x_{\wfc{\iota_n}{x},n}} - 1$ and
\[
\wfc{\rho_n}{x} = n \wlr{x - 2 \wlr{n \wlr{x + 1}/2 - \wfc{\theta_n}{x}}/n + 1}- 1
= 2 \wfc{\theta_n}{x} - 1 \in (-1,1).
\]
Therefore $\wfc{\rho^2_n}{x} < 1$,
and the definition of $D_n$ in Equation \pRef{numDen} leads to
\[
\wfc{D_n}{x} = \sum_{k = 0}^{\wfc{\iota_n}{x}} \wlr{-1}^k \frac{1}{x - x_{k,n}} +
               \sum_{k = \wfc{\iota_n}{x} + 1}^{n} \wlr{-1}^k \frac{1}{x - x_{k,n}}
\]
\[
 = \sum_{k = 0}^{\wfc{\iota_n}{x}}
 \wlr{-1}^{\wfc{\iota_n}{x} - k}
 \frac{1}{x - x_{\wfc{\iota_n}{x} - k,n}} +
 \sum_{k = 0}^{n - \wfc{\iota_n}{x} - 1} \wlr{-1}^{\wfc{\iota_n}{x} + k}
 \frac{1}{x_{\wfc{\iota_n}{x} + k + 1,n} - x}
\]
\[
 =
 \wlr{-1}^{\wfc{\iota_n}{x}}
 \wlr{ \sum_{k = 0}^{\wfc{\iota_n}{x}}
 \wlr{-1}^{k}
 \frac{1}{x - x_{\wfc{\iota_n}{x},n} + 2 k / n} +
 \sum_{k = 0}^{n - \wfc{\iota_n}{x} - 1} \wlr{-1}^{k}
 \frac{1}{x_{\wfc{\iota_n}{x} + 1,n} - x + 2 k / n}
 }
\]
\[
 =
 \wlr{-1}^{\wfc{\iota_n}{x}}
n \wlr{ \sum_{k = 0}^{\wfc{\iota_n}{x}}
 \wlr{-1}^{k}
 \frac{1}{ n \wlr{x - x_{\wfc{\iota_n}{x},n}} + 2 k} +
 \sum_{k = 0}^{n - \wfc{\iota_n}{x} - 1} \wlr{-1}^{k}
 \frac{1}{ n \wlr{x_{\wfc{\iota_n}{x},n} - x} + 2 k + 2}
 }
\]
\[
=
\wlr{-1}^{\wfc{\iota_n}{x}} n
 \wlr{ \sum_{k = 0}^{\wfc{\iota_n}{x}}
 \wlr{-1}^{k}
 \frac{1}{2 k + 1 + \wfc{\rho_n}{x}} +
 \sum_{k = 0}^{n - \wfc{\iota_n}{x} - 1} \wlr{-1}^{k}
 \frac{1}{2 k + 1 - \wfc{\rho_n}{x}}
 }.
\]
Therefore,
\pbDef{ddn}
\wfc{D_n}{x} = \wlr{-1}^{\wfc{\iota_n}{x}} n \wlr{ \wfc{U_n}{x} + \wfc{V_n}{x} }
\peDef{ddn}
for
\pbDef{bndn}
\wfc{U_n}{x} := \sum_{k = 0}^{\wfc{\iota_n}{x}}
 \wlr{-1}^{k}
 \frac{1}{2 k + 1 + \wfc{\rho_n}{x}}
\peDef{bndn}
and
\pbDef{cndn}
\wfc{V_n}{x} := \sum_{k = 0}^{n - \wfc{\iota_n}{x} - 1} \wlr{-1}^{k}
 \frac{1}{2 k + 1 - \wfc{\rho_n}{x}}.
\peDef{cndn}
Since $\wfc{\rho_n}{x} \in (-1,1)$ the absolute values of the parcels
of the sum $\wfc{U_n}{x}$ and $\wfc{V_n}{x}$ decrease with $k$,
their sign alternate, and the first parcel is positive.
Therefore, $\wfc{U_n}{x}$ and $\wfc{V_n}{x}$ are positive and Equation \pRef{ddn}
shows that $\wfc{D_n}{x}$ has the sign claimed by Lemma \ref{lemDen}.
Moreover, the definition \pRef{defA} of the function $A$ shows that
\[
\wfc{A}{x} = \sum_{k = 0}^\infty \wlr{-1}^{k}
\wlr{\frac{1}{2 k + 1 - \sqrt{x}} + \frac{1}{2 k + 1 + \sqrt{x}}}
\]
and Equation \pRef{ddn} yields
\[
\wabs{\wfc{D_n}{x}}/n - \wfc{A}{\wfc{\rho_n^2}{x}} =
\sum_{k = \wfc{\iota_n}{x} + 1}^{\infty}
 \wlr{-1}^{k}
 \frac{1}{2 k + 1 + \wfc{\rho_n}{x}} +
 \sum_{k = n - \wfc{\iota_n}{x}}^{\infty} \wlr{-1}^{k}
 \frac{1}{2 k + 1 - \wfc{\rho_n}{x}}.
\]
It follows that
\[
\wabs{\wabs{\wfc{D_n}{x}}/n - \wfc{A}{\wfc{\rho_n^2}{x}}}
\leq \wfc{G_n}{x} +  \wfc{H_n}{x},
\]
for
\[
\wfc{G_n}{x} :=
\sum_{k = 0}^{\infty}
 \wlr{-1}^{k}  \frac{1}{2 k + 2 \wfc{\iota_n}{x} + 3 + \wfc{\rho_n}{x}}
\]
and
\[
\wfc{H_n}{x} :=
\sum_{k = 0}^{\infty}
 \wlr{-1}^{k}  \frac{1}{2 k + 2 \wlr{n - \wfc{\iota_n}{x}} + 1 - \wfc{\rho_n}{x}}.
\]
Replacing $k$ by $2 \ell$ and $2 \ell+ 1$ in the expression of $G_n$ above
we obtain
\[
\wfc{G_n}{x} =
\sum_{\ell = 0}^{\infty}
 \frac{1}{4 \ell + 2 \wfc{\iota_n}{x} + 3 + \wfc{\rho_n}{x}} -
 \frac{1}{4 \ell + 2 \wfc{\iota_n}{x} + 5 + \wfc{\rho_n}{x}}
\]
\[
=
\sum_{k = 0}^{\infty}
 \frac{2}
 {
 \wlr{4 \ell + 2 \wfc{\iota_n}{x} + 3 + \wfc{\rho_n}{x}}
 \wlr{4 \ell + 2 \wfc{\iota_n}{x} + 5 + \wfc{\rho_n}{x}}
 }
\]
\[
\leq \int_{t = 0}^\infty
\frac{2}
 {
 \wlr{4 t + 2 \wfc{\iota_n}{x} + 3 + \wfc{\rho_n}{x}}
 \wlr{4 t + 2 \wfc{\iota_n}{x} + 5 + \wfc{\rho_n}{x}}
 }
 dt
\]
\[
= \frac{1}{4} \wfc{\ln}{1 + \frac{2}{2 \wfc{\iota_n}{x} + 3 + \wfc{\rho_n}{x}}}
\leq \frac{1}{2 \wlr{2 \wfc{\iota_n}{x} + 3 + \wfc{\rho_n}{x}}}
\leq \frac{1}{4 \wlr{1 + \wfc{\iota_n}{x}}}.
\]
The integral above was computed with Wolfram Alpha, and
a similar computation shows that
\[
\wfc{H_n}{x} \leq
\frac{1}{4 \wlr{n - \wfc{\iota_n}{x}}},
\]
and the second part of Equation \pRef{boundDen} holds.
It follows that
\[
\wabs{\wfc{D_n}{x}/n} \geq \wfc{A}{\wfc{\rho_n^2}{x}} - 1/2 
\geq \frac{\pi}{2} - 1/2 > 1.07 > 1,
\]
because $\wfc{A}{\wfc{\rho_{n_j}^2}{x}} \geq \pi/2$.
This proves the first part of bound \pRef{boundDenA}.
We also have 
\[
\wabs{\wfc{D_n}{x}/n}
\geq
\wfc{A}{\wfc{\rho_n^2}{x}} - 1/\pi \wfc{A}{\wfc{\rho_n^2}{x}}
> \wfc{A}{\wfc{\rho_n^2}{x}}/2,
\]
and $2/(1 - \wfc{\rho_n^2}{x}) \geq 2$ because
$\wfc{\rho_n^2}{x}\in [0,1)$. Equation
\pRef{boundA} shows that
\[
\wfc{A}{\wfc{\rho_n^2}{x}} \geq \frac{2}{1 - \wfc{\rho_n^2}{x}} - \frac{1}{2}
\geq \frac{2}{1 - \wfc{\rho_n^2}{x}} - \frac{1}{4} \times \frac{2}{1 - \wfc{\rho_n^2}{x}}
 = \frac{3}{4 \wlr{1 - \wfc{\rho_n^2}{x}}},
\]
and this proves the second Equation in \pRef{boundDenA}.

Finally, for every $x \in (-1,1)$ regular we have that
\[
\lim_{n_j \rightarrow \infty} \wfc{\iota_{n_j}}{x} =
\lim_{n_j \rightarrow \infty} n - \wfc{\iota_{n_j}}{x} =
 +\infty.
\]
This observation and the equations above imply
Equation \pRef{boundDenR}.
\peProof{Lemma}{lemDen}\\


\pbProofB{Corollary}{corConv}
Let us assume Equation \pRef{corConvLhs} and prove Equation \pRef{corConvRhs}.
Lemma \ref{lemDen} shows that $\wabs{\wfc{D_n}{x}/n} \geq 1$ for all
$x$ and $n$. Therefore, $L \neq 0$ and for $j$ large enough we must have
\[
\wfc{\wrm{sign}}{\wfc{D_{n_j}}{x}} = \wfc{\wrm{sign}}{L},
\]
and Equation \pRef{boundDen} shows that
this is also the sign of $\wlr{-1}^{\wfc{\iota_{n_j}}{x}}$. Therefore,
\[
\lim_{j \rightarrow \infty} \wlr{-1}^{\wfc{\iota_{n_j}}{x}} = \wfc{\wrm{sign}}{L}.
\]
Moreover, Equation \pRef{boundDenR} implies that
\pbTClaim{corConvB}
\lim_{j \rightarrow \infty} \wfc{A}{\wfc{\rho_{n_j}^2}{x}} = L.
\peTClaim{corConvB}
Since $A$ is continuous and $[0,1]$ is compact, this implies
that $L \in \wfc{A}{[0,1]} = [\pi/2,+\infty]$ and $L \geq \pi/2$.
Finally, since $A^{-1}$ is continuous Equation \pRef{corConvB}
implies that
$\lim_{j \rightarrow \infty} \wfc{\rho_n^2}{x} = \wfc{A^{-1}}{L}$
and the proof of Equation \pRef{corConvRhs} is complete.

Let us now assume Equation \pRef{corConvRhs} and prove
Equation \pRef{corConvLhs}. The continuity of
$A$ implies that
\[
\lim_{j \rightarrow \infty}
\wfc{A}{\wfc{\rho_{n_j}^2}{x}} = \wabs{L},
\]
and Equation \pRef{boundDenR} implies that
\[
\lim_{j \rightarrow \infty} \frac{1}{n_j} \wabs{\wfc{D_{n_j}}{x}} = \wabs{L}.
\]
This equation combined with the assumption
\[
\lim_{j \rightarrow \infty} \wlr{-1}^{\wfc{\iota_{n_j}}{x}} = \wfc{\wrm{sign}}{L}
\]
and Equation \pRef{boundDen} implies Equation \pRef{corConvLhs} and we are done.
\peProof{Corollary}{corConv}\\


\pbProofB{Lemma}{lemIrrat}
Let $r_i \in \wset{0,1}$ be such that $\wlr{-1}^{r_i} = \wfc{\wrm{sign}}{y}$
and
\[
z := \sqrt{\wfc{A^{-1}}{\wabs{y}}},
\]
so that if $z_j$ is a sequence such that
$\lim_{j \rightarrow \infty} z_j = z$ then
\[
y = \lim_{j \rightarrow \infty} \wlr{-1}^{r_i} \wfc{A}{ z_j^2 }.
\]
Lemma \ref{lemDen}
shows that to prove Lemma \ref{lemIrrat} it suffices
to define a sequence $n_j$ such that
\pbTClaim{reqIrrat}
n_j \equiv \ r_n \ \wrm{mod} \ 2,
\hspace{1cm}
\wfc{\iota_{n_j}}{x} \equiv r_i \ \wrm{mod} \ 2
\hspace{1cm} \wrm{and} \hspace{0.5cm}
\lim_{j \rightarrow \infty} \wfc{\rho_{n_j}}{x} = z.
\peTClaim{reqIrrat}
Since the image of $A^{-1}$ is $[0,1]$ we have that $z \in [0,1]$,
and there exist sequences $p_j,q_j \in \wn{}$
with $\lim_{j \rightarrow \infty} p_j / q_j = z$
and $0 < p_j/q_j < 1$. We start with an empty set of integers $n_j$,
and build them by induction. At the $j$th step
we use Hartmann's Theorem with
$\xi = x + 1$, $s = 4 q_j$, $a = 2 q_j r_i + p_j + q_j$
and $b = q_j r_n$,
and conclude that there exist infinitely many numbers $u$ and $v$ such
that
\[
\wabs{x + 1 - \frac{4 q_j u + 2 q_j r_i + p_j + q_j}{4 q_j  v + q_j r_n}} <
\frac{32 q_j^2}{\wlr{4 q_j v + q_j r_n}^2} \leq \frac{2}{v^2}.
\]
This implies that
\pbTClaim{xIrrat}
x + 1 = \frac{4 u + 2 r_i + 1 + p_j/q_j}{4 v +  r_n}
+ \theta_{j} \frac{1}{v^2}
\peTClaim{xIrrat}
for some $\theta_{j} \in [-2,2]$.
Taking a pair $\wlr{u_j,v_j}$ with $v_j$ so large that
\pbTClaim{largeV}
0 < 1 + p_j/q_j - 2 \wlr{4 v_j + r_n}/v_j^2 < 1 + p_j/q_j + 2 \wlr{4 v_j + r_n}/{v_j^2} < 2,
\peTClaim{largeV}
and for which $n_j := 4 v_j + r_n$ is larger than the previous $n_j$, we obtain a
$n_j$ which satisfies the parity requirement in Lemma \ref{lemIrrat} and
\[
n_j \wlr{x + 1}/2 = 2 u_j + r_i + \wlr{1 + p_j/q_j + \theta_{j} n_j/v_j^2}/2.
\]
The definition \pRef{rho} of $\iota_n$ and Equation \pRef{largeV} implies that
\[
\wfc{\iota_{n_j}}{x} = \wfloor{n_j \wlr{x + 1} / 2} = 2 u_j + r_i,
\]
and this $\wfc{\iota_{n_j}}{x}$ has the parity claimed by Equation \pRef{reqIrrat}, and
Equation \pRef{xIrrat} yields
\[
x - x_{\wfc{\iota_{n_j}}{x},n_j} = x + 1 - 2 \frac{2 u_j + r_i}{4 v_j + r_n}
= \frac{1 + p_j/q_j}{4 v_j + r_n} + \theta_{j} \frac{1}{v_j^2},
\]
and the definition \pRef{rho} of $\rho$ yields
\[
\wfc{\rho_{n_j}}{x} = n_j \wlr{x - x_{\wfc{\iota_{n_j}}{x},n_j}} - 1
= p_j/q_j + \theta_j \wlr{4 v_j + r_n}/v_j^2.
\]
Since $\wabs{\theta_j} \leq 2$ and $r_n \in \wset{0,1}$,
we have that
\[
\lim_{j \rightarrow \infty} \wfc{\rho_{n_j}}{x} =
\lim_{j \rightarrow \infty} p_j/q_j = z,
\]
and the proof of Lemma \ref{lemIrrat} is complete.
\peProof{Lemma}{lemIrrat}\\

%
%

\pbProofB{Lemma}{lemRat}
Corollary \ref{corConv} shows that
\pbTClaim{rhoM}
\lim_{j \rightarrow \infty} \wabs{\wfc{\rho_{n_j}}{x}} = M := \wfc{A^{-1}}{\wabs{L}},
\peTClaim{rhoM}
and there exists $j_0$ such that if $j \geq j_0$ then
$\wlr{-1}^{\wfc{\iota_{n_j}}{x}} = \wfc{\wrm{sign}}{L}$ and
\pbDef{epsRat}
\wabs{\wfc{\rho_{n_j}}{x}} = M + \epsilon_j
\hspace{1cm} \wrm{with} \hspace{1cm}
\wabs{\epsilon_j} \leq \frac{1 - \wlr{q M - \wfloor{q M}}}{2 q}.
\peDef{epsRat}

Equation \pRef{rhoDec} and the hypothesis $x = p/q - 1$ imply that
\pbTClaim{iotaM}
p/q = \wlr{2 \wfc{\iota_{n_j}}{x} + \sigma_j \wlr{M + \epsilon_j} + 1}/n_j,
\peTClaim{iotaM}
with $\sigma_j \in \wset{-1,1}$, and
\pbTClaim{epsInteger}
p n_j - 2 q \wfc{\iota_{n_j}}{x} - \sigma_j \wfloor{q M} - q =
\sigma_j \wlr{q \epsilon_j + \wlr{q M - \wfloor{q M}}}.
\peTClaim{epsInteger}
Since $\wabs{\sigma_j} = 1$,
Equation \pRef{epsRat} yields
\[
\wabs{\sigma_j  \wlr{q \epsilon_j + \wlr{q M - \wfloor{q M}}}}
\leq \wlr{1 - \wlr{q M - \wfloor{q M}}}/2 + q M - \wfloor{q M}
\]
\[
= \wlr{1 + \wlr{q M - \wfloor{q M}}}/2 < 1.
\]
Since the left hand side of Equation \pRef{epsInteger} is
integer, we have that
\[
q \epsilon_j + \wlr{q M - \wfloor{q M}} = 0 \Rightarrow
\epsilon_j = \wfloor{q M}/q - M,
\]
Equation \pRef{epsRat} yields
\[
\wabs{\wfc{\rho_{n_j}}{x}} = \wfloor{q M}/q,
\]
and Equation \pRef{rhoM} implies that
$\wfloor{q M} = q M$. It follows that
$q M \in \wz{}$ and $M = m/q$ for some $m \in \wz{}$.
Therefore,
$\wabs{\wfc{\rho_{n_j}}{x}} = m/q$,
and the proof of Lemma \ref{lemRat} is complete.
\peProof{Lemma}{lemRat}\\

%
%

\pbProofB{Corollary}{corRatPQ}
For a regular $x = p/q - 1$,
with $\lim_{j \rightarrow \infty} \wfc{D_{n_j}}{x}/n_j = L$,
Lemma \ref{lemRat} implies that there exist
$i_j \in \wn{}$,  $m \in \wz{}$ with $\wabs{m} < q$, and $s \in \wset{0,1}$ such that
\[
p/q = \wlr{2 \wlr{2 i_j + s} + m/q + 1}/n_j,
\hspace{0.5cm}
\wlr{-1}^s = \wfc{\wrm{sign}}{L}
\hspace{0.5cm} \wrm{and} \hspace{0.5cm}
\wabs{L} = \wfc{A}{m^2/q^2},
\]
and the first Equation above is equivalent to
\pbDef{ratPqA}
p \, n_j = 2 \wlr{2 i_j + s} q + m + q.
\peDef{ratPqA}
When $n_j$ is odd, this equation implies that
$\ell :=  m/2 \in \wz{}$, and $\wabs{\ell} \leq \wlr{q - 1}/2$.
Therefore, $\wabs{L} = \wfc{A}{4 \ell^2 / q^2}$
and the set in Equation \pRef{oddRatPQ} does contain all
relevant limits $L$.
Conversely, with $m = 2\ell$ and $m_j = \wlr{n_j - 1}/2$, Equation \pRef{ratPqA}
is equivalent to
\[
p \, m_j =  \wlr{2 i_j + s} q + \ell + \wlr{q - p}/2
= \wlr{2 q} i_j + \wlr{s q + \ell + \wlr{q - p}/2}.
\]
For every $s$ and $\ell$ this equation has
infinitely many solutions $\wlr{m_j,i_j} \in \wn{} \times \wn{}$ because
$\wfc{\wrm{gcd}}{p,2q} = 1$. Therefore, for every $m = 2 \ell$,
and $s \in \wset{0,1}$ there exist infinitely many $n_j = 2 m_j + 1$
which satisfy Equation \pRef{ratPqA}, and all elements
in the set $\wfc{\wcal{O}}{p/q}$ in Equation \pRef{oddRatPQ}
are indeed limits of sequences $\wfc{D_{n_j}}{x}/n_{j}$
with odd $n_j$.
This completes the verification of Equation \pRef{oddRatPQ}.

When $n_j$ is even, Equation
\pRef{ratPqA} implies that $\ell := \wlr{m - 1}/2 \in \wz{}$,
$\wabs{2 \ell + 1} < q$
and $\wabs{L} = \wfc{A}{\wlr{2 \ell + 1}^2 / q^2}$,
and the set in Equation \pRef{evenRatPQ} does contain
all the relevant limits $L$.
Moreover, for $m_j = n_j/2 \in \wz{}$ and $m = 2 \ell + 1$,
Equation \pRef{ratPqA} reduces to
\[
p \,  m_j =  \wlr{2 i_j + s} q + \ell + \wlr{q + 1}/2 = \wlr{2 q} i_j +
\wlr{s q + \ell + \wlr{q + 1}/2}
\]
and, as before, we can find infinitely many $\wlr{m_j,i_j}$ which
satisfy this equation, and use then to generate sequences $n_j$
with all the limits in the set in Equation \pRef{evenRatPQ}.
As a result, Equation \pRef{evenRatPQ} is valid, and this proof is complete.
\peProof{Corollary}{corRatPQ}\\


\pbProofB{Corollary}{corRat2PQ}
If $x = 2 p/q - 1$ is regular and
$\lim_{j \rightarrow \infty} \wfc{D_{n_j}}{x}/n_j = L$ then
Lemma \ref{lemRat} implies that there exist
$i_j \in \wn{}$, $m \in \wz{}$ with $\wabs{m} < q$ and $s \in \wset{0,1}$ such that
\[
2p/q = \wlr{2 \wlr{2 i_j + s} + m/q + 1}/n_j,
\hspace{0.5cm}
\wlr{-1}^s = \wfc{\wrm{sign}}{L}
\hspace{0.5cm} \wrm{and} \hspace{0.5cm}
\wabs{L} = \wfc{A}{m^2/q^2}.
\]
The first Equation above is equivalent to
\[
2 p n_j = 2 \wlr{2 i_j + s} q + m + q,
\]
and it implies that $h := \wlr{m - 1}/2 \in \wz{}$.
Therefore,
\pbDef{rat2PqA}
p \, n_j = \wlr{2 i_j + s} q + h + \wlr{q+1}/2.
\peDef{rat2PqA}

If $n_j$ is odd then $\ell := \wlr{s + h - p + \wlr{q+1}/2} \in \wz{}$,
and $m = 4 \ell + 2 p - 2 s - q$.
Since $\wabs{m} < q$ we have that
\[
-q + 1 \leq 4 \ell + 2 p - 2 s - q \leq q - 1
\]
and
\[
1 - 2 p + 2 s \leq 4 \ell \leq - 2 p + 2 s + 2 q - 1
\ \Rightarrow \
s - p + 1 \leq 2 \ell  \leq  s - p + q - 1,
\]
and the set in Equation \pRef{oddRat2PQ} contains
all the relevant limits.
Conversely, for $m_j := \wlr{n_j + 1}/2$ and $h = 2 \ell + p - s - \wlr{q + 1}/2$,
 Equation \pRef{rat2PqA} reduces to
\[
p m_j = i_j q + \ell + s \wlr{q - 1}/2,
\]
and since $\wfc{\wrm{gcd}}{p,q} = 1$ there exist
infinitely many $m_j$ and $i_j$ which satisfy
this equation, and all elements of the
set $\wfc{\wcal{O}}{2p/q}$ in Equation \pRef{oddRat2PQ} are indeed
limits corresponding to conveniently chosen odd sequences.

If $n_j$ is even then Equation \pRef{rat2PqA} yields
$\ell := \wlr{s + h + \wlr{q+1}/2} \in \wz{}$.
Since $m = 2 h + 1$, we obtain
\[
h = 2 \ell - s - \wlr{q+1}/2
\ \ \Rightarrow \ \
m = 4 \ell - 2 s - q,
\]
and the bound $\wabs{m} < q$ leads to
$1 +  s  \leq 2 \ell < s + q - 1$,
and Equation \pRef{evenRat2PQ} is correct.
Finally, with $m_j := n_j /2 \in \wz{}$ and $h$ above, Equation \pRef{rat2PqA}
reduces to
\[
p m_j = i_j q + \ell + s \wlr{q - 1}/2,
\]
and since $\wfc{\wrm{gcd}}{p,q} = 1$ there exist
infinitely many $\wlr{m_j,i_j}$ which satisfy
this equation.
\peProof{Corollary}{corRat2PQ}\\


\pbProofB{Corollary}{corRatP2Q}
If $x = p/{2q} - 1$ is regular and
$\lim_{j \rightarrow \infty} \wfc{D_{n_j}}{x}/n_j = L$ then
Lemma \ref{lemRat} implies that there exist
$i_j \in \wn{}$, $m \in \wz{}$ with $\wabs{m} < 2 q$, and $s \in \wset{0,1}$, such that
\[
\frac{p}{2 q} = \wlr{2 \wlr{2 i_j + s} + m/{2q} + 1}/n_j,
\hspace{0.5cm}
\wlr{-1}^s = \wfc{\wrm{sign}}{L}
\hspace{0.5cm} \wrm{and} \hspace{0.5cm}
\wabs{L} = \wfc{A}{\frac{m^2}{4 q^2}}.
\]
The first equation above is equivalent to
\pbDef{ratP2qA}
p \, n_j = 4 \wlr{2 i_j + s} q + m + 2 q.
\peDef{ratP2qA}
When $n_j$ is odd, $\ell := \wlr{m - 1} / 2 \in \wz{}$ and
the bound $\wabs{m} < 2 q$ implies that $\wabs{\ell} \leq q - 1$
and Equation \pRef{oddRatP2Q} is correct.
Conversely, for $m = 2 \ell + 1$ and $m_j = \wlr{n_j - 1}/2$
Equation \pRef{ratP2qA} reduces to
\[
p \, m_j +  =\wlr{4 q} i_j + s q + \ell + \frac{1 - p}{2} + q,
\]
and since $\wfc{\wrm{gcd}}{p,4 q} = 1$,
for each $s$ and $\ell$
this equation has infinitely many solutions
$(n_j,i_j)$, which we can use to build sequences
with the limits in the set in Equation \pRef{oddRatP2Q}.

When $n_j$ is even, Equation \pRef{ratP2qA}
implies that $\ell := m / 2 \in \wz{}$ and
the bound $\wabs{m} < 2 q$ implies that $\wabs{\ell} \leq q - 1$
and Equation \pRef{evenRatP2Q} is correct.
Conversely, for $m = 2 \ell$ above and $m_j := n_j/2 \in \wz{}$,
Equation \pRef{ratP2qA} reduces to
\[
p \, m_j = 4 i_j + 2 s q + \ell + q,
\]
and since $\wfc{\wrm{gcd}}{p,4q} = 1$,
for each $s$ and $\ell$
this equation has infinitely many solutions
$(n_j,i_j)$, from which we can obtain sequences
with the limits in Equation
\pRef{evenRatP2Q}.

\peProof{Corollary}{corRatP2Q}\\


\section{The numerator of the error for $f$ in $\wrm{AC}^1$}
\label{secNum}
In this section we explore the consequences
of the observation in the introduction
that Berrut's interpolants are biased.
After we remove the bias, the relevant quantity
for understanding the convergence of the interpolants
$B_n$ is defined as
\pbDef{diffN}
\wfc{\Delta_{n}}{f,x} :=
\frac{\wfc{f}{-1} - \wfc{f}{x}}{2 \wlr{x + 1}} +
\wlr{-1}^n \frac{\wfc{f}{1} - \wfc{f}{x}}{2 \wlr{x - 1}} +
\frac{1}{n} \sum_{k = 1}^{n-1} \wlr{-1}^k \frac{\wfc{f}{x_{k,n}} - \wfc{f}{x}}{x - x_{k,n}}
\peDef{diffN}
for $x \not \in \wset{x_{0,n},\dots,x_{nn}}$, and
$\wfc{\Delta_{n}}{f,x_{k,n}} := 0$.
We can then express the combination of $\wfc{\Delta_n}{f,x}$ and the bias
$\wfc{O}{f,x}$ for $n_j = 2 j + 1$ odd as
\pbDef{diffNO}
\wfc{B_{2n+1}}{f,x} - \wfc{f}{x} =
\wlr{\wfc{\Delta_{2n+1}}{f,x} + \wfc{O}{f,z}} / \wfc{D_{2n+1}}{x}.
\peDef{diffNO}
For $n_j = 2 n$ the bias is $\wfc{E}{f,x}$ and we have
\pbDef{diffNE}
\wfc{B_{2n}}{f,x} - \wfc{f}{x} =
\wlr{\wfc{\Delta_{2n}}{f,x} + \wfc{E}{f,x}} / \wfc{D_{2n}}{x}.
\peDef{diffNE}
The expression for $\wfc{\Delta_n}{f,x}$ for both parities
is the same, that is, the bias is related to parity,
but the mean term $\wfc{\Delta_n}{f,x}$ is not.
We can then obtain a clean result regarding the convergence
of the numerator of the error, which we prove in the
end of this section.

\pbTheoremBT{thm2}{The uniform convergence of the numerator}
If $f \in \wrm{AC}^1$ then
\[
\lim_{n \rightarrow \infty} \,
\wnorm{\wfc{\Delta_{n}}{f}}_\infty = 0.
\]
\peTheorem{thm2}


\pbProofB{Theorem}{thm2}
Given $\epsilon \in (0,1)$, by the absolute continuity of $\wdfc{f}{x}$ there exists $\delta \in (0,1)$
for which
\pbDef{acdelta}
\sum_{k = 0}^m \wabs{y_k  - z_k} \leq \delta
 \Rightarrow \sum_{k=0}^m \wabs{\wdfc{f}{y_k} - \wdfc{f}{z_k}} < \epsilon/3,
\peDef{acdelta}
and we now define
\pbDef{defn0}
n_0 :=
1024 + \wceilB{\wlr{1 +
\wnorm{f''}_1^2} \frac{1}{100 \, \delta \,  \epsilon^2}},
\peDef{defn0}
take $n \geq n_0$ and $x \in [-1,1]$
and show that $\wabs{ \wfc{\Delta_{n}}{f,x} } \leq \epsilon$.
If $x \in \wset{x_{0n},\dots,x_{n,n}}$ then
$\wfc{\Delta_{n}}{f,x} = 0$ by definition and we are done.
For $x \not \in \wset{x_{0n},\dots,x_{n,n}}$,
let $i$ be the index such that $x_{i,n}$ is the node closest to $x$.
We split the sum which defines $\wfc{\Delta_n}{f,x}$
in Equation \pRef{diffN} in at most
three parts: $F$ (first), $M$ (middle) and $L$ (last),
according to the distance of $x$ to $\pm 1$.
When $x$ is too close to $-1$ we leave the First region
empty, and if $x$ is too close to $1$ then the Last range is
left empty. When not empty, the First range corresponds
to parcels with indexes from $0$ to $2m$.
The Middle range spans the indexes from $2m$ to $n - 2\ell$,
and contains of the order of $\sqrt{\delta n}$ parcels
(the parcel corresponding to $k = 2m$ is split between
the First and Middle ranges.) When not empty, the Last
range starts at index $n - 2\ell$ and ends a index $n$,
and the parcel of index $n - 2 \ell$ is split between
the Middle and Last ranges.

Formally, we define
\begin{enumerate}
\item If $i < \sqrt{\delta n}/4$ then, since $n \geq n_0$,
Equation \pRef{defn0} implies that $n > i + \sqrt{\delta n}/4$ and we define
$m := 0$,
\pbDef{defLI}
\ell = n - 2 \wfloorB{\wlr{n - i - \sqrt{\delta n}/4}/2} + 2,
\peDef{defLI}
$F := 0$,
\[
M :=
\frac{\wfc{f}{x_{2m,n}} - \wfc{f}{x}}{2 \wlr{x - x_{2m,n}}} + \wlr{-1}^{n}
\frac{\wfc{f}{x_{n - 2 \ell,n}} - \wfc{f}{x}}{2 \wlr{x - x_{n - 2 \ell,n}}}
\]
\pbDef{defM}
+ \sum_{k = 2m + 1}^{n - 2 \ell - 1}
\wlr{-1}^k \frac{\wfc{f}{x_{k,n}} - \wfc{f}{x}}{x - x_{k,n}}
\peDef{defM}
and
\[
L :=
\wlr{-1}^n  \frac{\wfc{f}{x_{n - 2 \ell,n}} - \wfc{f}{x}}{2 \wlr{x - x_{n - 2 \ell,n}}} +
\wlr{-1}^n   \frac{\wfc{f}{1} - \wfc{f}{x}}{2 \wlr{x - 1}} +
\]
\pbDef{defL}
\sum_{k = n - 2 \ell + 1}^{n} \wlr{-1}^k \frac{\wfc{f}{x_{k,n}} - \wfc{f}{x}}{x - x_{k,n}}.
\peDef{defL}
\item If $\sqrt{\delta n}/4 \leq i < n - \sqrt{\delta n}/4$ then
we set
\pbDef{defMid}
m  := 2 \wceilB{\wlr{i - \sqrt{\delta n}/4}/2} + 2,
\peDef{defMid}
define $M$ and $L$ as in Equations \pRef{defL} and \pRef{defMid}, and
\pbDef{defF}
F :=
\frac{\wfc{f}{-1} - \wfc{f}{x}}{2 \wlr{x + 1}} +
\frac{\wfc{f}{x_{2 m,n}} - \wfc{f}{x}}{2 \wlr{x - x_{2m,n}}}
+ \sum_{k = 1}^{2 m - 1}
\wlr{-1}^k \frac{\wfc{f}{x_{k,n}} - \wfc{f}{x}}{x - x_{k,n}}.
\peDef{defF}
\item Finally, if $i \geq n - \sqrt{\delta n}/4$ then we
define $m$ as in Equation \pRef{defMid}, $\ell = 0$,
 $M$ and $F$ as in Equations
\pRef{defL}  and \pRef{defF}, and $L := 0$.
\end{enumerate}
We now bound $M$. Splitting each parcel in two parts,
and grouping consecutive halves
and using the Mean Value Theorem we obtain
\[
2 \wabs{M} =
\wabs{\sum_{k = 2 m}^{n - 2 \ell - 1}
\wlr{\frac{\wfc{f}{x_{k,n}} - \wfc{f}{x}}{x - x_{k,n}} -
\frac{\wfc{f}{x_{k+1,n}} - \wfc{f}{x}}{x - x_{k+1,n}}}}
\leq
\sum_{k = 2m}^{n - 2 \ell - 1} \wabs{\wdfc{f}{\xi_k} - \wdfc{f}{\xi_{k+1}}}
\]
with
\[
\wabs{\xi_k - x} \leq 2 \max \wset{i - 2 m, n - 2\ell - i}/n.
\]
The indexes $\ell$ and $m$ were defined in
Equations \pRef{defLI} and \pRef{defMid} so that
\[
0 < i - 2m \leq \sqrt{\delta n}/4,
\hspace{0.5cm}
0 < n - 2\ell - i \leq \sqrt{\delta n}/4
\hspace{1cm} \wrm{and} \hspace{1cm}
\wabs{\xi_k - x} \leq \frac{\sqrt{\delta n}}{2n}
\]
This implies that $\wabs{\xi_k - \xi_{k+1}} \leq \sqrt{\delta/ n}$,
\[
\sum_{k = 2m}^{n - 2 \ell} \wabs{\xi_k - \xi_{k+1}}
 \leq \wlr{n - 2 \ell - 2 m} \sqrt{\delta/n}
 \leq \sqrt{\delta n}/2 \times  \sqrt{\delta/n}
 \leq \delta,
\]
and Equation \pRef{acdelta} implies that
$M \leq \epsilon / 3$.
We now show that
$L \leq \epsilon/3$ in the case in
which it is different from zero (By symmetry,
the same bound applies to $F$.)

Defining $y_k = x_{n - 2 \ell + k,n}$,
we can group the terms of $L$ as
\[
-2 \wlr{-1}^n L = \sum_{j = 0}^{\ell - 1}
\left( \ \
\wlr{
  \frac{\wfc{f}{y_{2j + 1}} - \wfc{f}{x}}{y_{2j + 1} - x} -
  \frac{\wfc{f}{y_{2j + 2}} - \wfc{f}{x}}{y_{2j + 2} - x}
  }
  \right.
\]
\[
\left.
  -
  \wlr{
  \frac{\wfc{f}{y_{2j}}     - \wfc{f}{x}}{y_{2 j} - x} -
  \frac{\wfc{f}{y_{2j + 1}} - \wfc{f}{x}}{y_{2 j + 1} - x}
  }
\ \ \right)
\]
\[
=
\frac{2}{n} \sum_{j = 0}^{\ell - 1}
[y_{2 j}, \, x, \, y_{2 j +1}, \, f] - [y_{2j + 2}, \, x, \, y_{2j +1}, \, f],
\]
where $[x_1, \, x_2,  \,  x_3,  \, f]$ denotes
the divided difference of second order corresponding to
$x_1$, $x_2$, $x_3$ and $f$, because
\[
y_{2j} - y_{2j + 1} = y_{2j + 1} - y_{2j + 2} = -2/n.
\]

Since $f'$ is absolutely continuous,
the Genocchi-Hermite formula \cite{deBoor} yields
\[
[u, \, v, \,  w,  \,  f] =
\int_0^1 \int_0^{1 - t} \wdsf{f}{\wlr{1 - t - s} u + s \, v + t \,  w} ds \ dt,
\]
and leads to
\pbDef{defLH}
-2 \wlr{-1}^n L = \frac{1}{n} \sum_{j = 0}^{\ell  - 1} \int_0^1 \wfc{h_j}{t} dt,
\peDef{defLH}
with $z_j := y_j - x > 0$ and
\[
\wfc{h_{j}}{t} :=
\int_0^{1 - t}
\wdsf{f}{x + t \, z_{2j} +  s \, z_{2j + 1}} ds
-
\int_0^{1 - t}
\wdsf{f}{x + t \, z_{2j + 2} + s \, z_{2j + 1}} ds.
\]
The changes of variables
\[
u = x + t \, z_{2j} +  s \, z_{2j + 1}
\hspace{0.5cm} \wrm{and} \hspace{0.5cm}
v = x + t \, z_{2j + 2} + s \,  z_{2j + 1}
\]
have the same Jacobian $z_{2 j + 1}$ with respect to $s$ and
\[
z_{2 j + 1} \, \wfc{h_{j}}{t} =
\int_{
x + t \, z_{2 j}
}^{
x + t \, z_{2j} + \wlr{1 - t} z_{2j + 1}
}
\wdsf{f}{u} du
-
\int_{
x + \,  t z_{2j + 2}
}^{
x +  \, t z_{2j + 2} + \wlr{1 - t } z_{2j + 1}
}
\wdsf{f}{v} dv.
\]
Since $z_{k + 1} - z_{k} = 2 / n$,
\[
z_{2 j + 1} \, \wfc{h_{j}}{t}  =
\int_{
x + t \, z_{2j}
}^{
x + z_{2j + 1} - 2 t/n
}
\wdsf{f}{u} du
-
\int_{
x + t \, z_{2j + 2}
}^{
x + z_{2j + 1} + 2 t/ n
}
\wdsf{f}{u} du
\]
\[
=
\int_{
x + t \, z_{2j}
}^{
x + t \, z_{2j + 2}
}
\wdsf{f}{u} du
-
\int_{
x + z_{2j + 1} - 2 t/n
}^{
x + z_{2j + 1} + 2t/ n
}
\wdsf{f}{u} du.
\]
The bound
\[
z_{2j + 1} = x_{2 j + n - 2 \ell,n} - x
\geq x_{n - 2 \ell,n} - x \geq \frac{2 \sqrt{\delta n} - 4}{4 n} \geq \sqrt{\delta/n}/3
\]
and the fact that $t \in [0,1]$ lead to
\[
\wabs{\wfc{h_{j}}{t}}
\leq
3 \sqrt{n/\delta}
\wlr{
\int_{
x + t \, z_{2j}
}^{
x + t \, z_{2j + 2}
}
\wabs{\wdsf{f}{u}} du
+
\int_{
x + z_{2j}
}^{
x + z_{2j + 2}
}
\wabs{\wdsf{f}{u}} du
}.
\]
It follows that
\[
\sum_{j = 0}^{\ell-1} \wabs{\wfc{h_j}{t}}
\leq  3 \sqrt{n/\delta}
\wlr{
\int_{
x + t \, z_{0}
}^{
x + t \, z_{2\ell}
}
\wabs{\wdsf{f}{u}} du
+
\int_{
x + z_{0}
}^{
x + z_{2\ell}
}
\wabs{\wdsf{f}{u}} du
}
\leq 6 \sqrt{n/\delta} \wnorm{f''}_1,
\]
and the same bound applies to $\int_0^1 \sum_{j = 0}^{\ell - 1} \wabs{\wfc{h_j}{t}} dt$.
The choice of $n_0$ in Equation
\pRef{defn0} and Equation \pRef{defLH} yield
\[
\wabs{L} \leq  3 \sqrt{\frac{1}{n \delta}} \wnorm{f''}_1
\leq \epsilon / 3,
\]
and we are done.
\peProof{Theorem}{thm2}\\

\section{The numerator of the error for $f$ in $BV^1$}
\label{secBV}
In this section we analyze the function $\wfc{\Delta_n}{f,x}$ 
defined in Equation \pRef{diffN} 
for functions $f$ with derivatives of bounded variation.
In summary, we show that in this case $\Delta_n$ is bounded
by half the total variation of $f'$. 
Our proof follows from this version of the Mean Value Theorem:

\pbTheoremBT{thmBV}{A monotone Mean Value Theorem}
Let $a, b \in \wrone{}$ be such that $a < b$, and
let $f: [a,b] \rightarrow \wrone{}$ be a continuous function,
which is differentiable in $(a,b)$. If $c$ and $\xi_c$ are such that
$a < \xi_c < c < b$ and
\pbDef{bvxic}
\wdfc{f}{\xi_c} = \frac{\wfc{f}{c} - \wfc{f}{a}}{c - a}
\peDef{bvxic}
then there exists $\xi_b \in [\xi_c,b)$ such that
\pbDef{bvxib}
\wdfc{f}{\xi_b} = \frac{\wfc{f}{b} - \wfc{f}{a}}{b - a}.
\peDef{bvxib}
\peTheorem{thmBV}
We prove Theorem \ref{thmBV} at the end of this section.
By induction, we conclude from this theorem that given an
increasing sequence $b_0,\dots,b_m$, with $b_0 > a$, 
we can find a non decreasing sequence $\xi_i$, with $\xi_i \in (a,b_i)$,
such that
\[
\wdfc{f}{\xi_i} = \frac{\wfc{f}{b_i} - \wfc{f}{a}}{b_i - a}.
\]
Using this observation, it is easy to prove the following
corollary:

\pbCorolBT{corNumBV}{The numerator of the error for $f$ in $BV^1$}
If $f \in \wrm{BV}^1$ and its derivative
has total variation $T_{f'}[-1,1] < +\infty$
then the function $\Delta_n$ in Equation \pRef{diffN} satisfies
\[
\wabs{\wfc{\Delta_n}{f,x}} \leq T_{f'}[-1,1] / 2.
\]
\peCorol{corNumBV}

In fact, if $x \in \wset{x_0,\dots,x_n}$, 
then $\wfc{\Delta_n}{f,x} = 0$ by definition.
For $x \not \in \wset{x_{0,n},\dots,x_{nn}}$, 
Equation \pRef{diffN} leads to 
\[
\wfc{\Delta_{n}}{f,x} =
\frac{1}{2} \sum_{k = 0}^{n-1} 
\wlr{
\frac{\wfc{f}{x_{k,n}} - \wfc{f}{x}}{x - x_{k,n}}
- 
\frac{\wfc{f}{x_{k + 1,n}} - \wfc{f}{x}}{x - x_{k+1,n}}
},
\]
and Theorem \ref{thmBV} yields an increasing sequence $\xi_0, \dots, \xi_n
\in [-1,1]$ such that
\[
\wdfc{f}{\xi_k} = -\frac{\wfc{f}{x_{k,n}} - \wfc{f}{x}}{x - x_{k,n}}.
\]
It then follows that
\[
\wabs{\wfc{\Delta_{n}}{f,x}} \leq
\frac{1}{2} \sum_{k = 0}^{n-1}
\wabs{
\frac{\wfc{f}{x_{k,n}} - \wfc{f}{x}}{x - x_{k,n}}
-
\frac{\wfc{f}{x_{k + 1,n}} - \wfc{f}{x}}{x - x_{k+1,n}}
}
\]
\[
= \frac{1}{2} \sum_{k = 0}^{n-1}
\wabs{\wdfc{f}{\xi_k} - \wdfc{f}{\xi_{k+1}}}
\leq  T_{f'}[-1,1] / 2,
\]
This proves Corollary \ref{corNumBV},
and we now present the proof of Theorem \ref{thmBV}.

\pbProofB{Theorem}{thmBV}
Let us start the proof with the particular case in which
\pbDef{bvc}
\wfc{f}{c} = \wfc{f}{a}.
\peDef{bvc}
By the traditional Mean Value Theorem, there exists $\mu \in \wlr{c, b}$ such that
\pbDef{bvu}
\wfc{f}{b} - \wfc{f}{c} = \wdfc{f}{\mu} \wlr{b - c} = v \wlr{b - a}
\hspace{1cm} \wrm{for} \hspace{1cm}
v := \wdfc{f}{\mu} \frac{b - c }{b - a}.
\peDef{bvu}
Equations \pRef{bvxic} and \pRef{bvc} imply that $\wdfc{f}{\xi_c} = 0$, and
since
\[
0 < \frac{b - c}{b - a} < 1
\]
we have that $v$ lies between $0 = \wdfc{f}{\xi_c}$
and $\wdfc{f}{\mu}$. Since derivatives have the intermediate value property,
there exists $\xi_b \in [\xi_c,\mu] \subset [\xi,b)$ such that
$\wdfc{f}{\xi_b} = v$. As a result, Equations \pRef{bvc} and \pRef{bvu} lead to
\[
\wdfc{f}{\xi_b} \wlr{b - a} = \wfc{f}{b} -\wfc{f}{c} = \wfc{f}{b} - \wfc{f}{a},
\]
and we have verified Equation \pRef{bvxib}
assuming that \pRef{bvc} holds.
To handle the general case it suffices to apply the argument above to
\pbDef{bvg}
\wfc{g}{x} = \wfc{f}{x} - \wlr{x - a} \frac{\wfc{f}{c} - \wfc{f}{a}}{c - a}.
\peDef{bvg}
In fact, $\wfc{g}{c} = \wfc{f}{a}= \wfc{g}{a}$ and Equation \pRef{bvxic} implies that
\[
\wdfc{g}{\xi_c} = \wdfc{f}{\xi_c} - \frac{\wfc{f}{c} - \wfc{f}{a}}{c - a} = 0 =
\frac{\wfc{f}{a} - \wfc{f}{a}}{c - a} =
\frac{\wfc{g}{c} - \wfc{g}{a}}{c - a}.
\]
As a result, the argument above yields $\xi_b \in [\xi_c,b)$ such that
\[
\frac{\wfc{g}{b} - \wfc{g}{a}}{b - a} =
\wdfc{g}{\xi_b} = \wdfc{f}{\xi_b} - \frac{\wfc{f}{c} - \wfc{f}{a}}{c - a}.
\]
It then follows from Equation \pRef{bvg} that
\[
\wdfc{f}{\xi_b} =
\frac{1}{b - a} \wlr{\wfc{g}{b} - \wfc{g}{a} +
     \wlr{b - a} \frac{\wfc{f}{c} - \wfc{f}{a}}{c - a}}
= \frac{\wfc{f}{b} - \wfc{f}{a}}{b - a},
\]
and we are done with the general case.
\peProof{Theorem}{thmBV}\\

\section{Combining the numerator with the denominator}
\label{secSecond}

In this section we combine the results from the previous
sections to prove Theorems \ref{thmMain}, \ref{thmUnif} and \ref{thmUnifBV}.


\pbProofB{Theorem}{thmMain}
We start with an odd sequence $n_j$ and an irrational $x$
for which
$\lim_{j \rightarrow \infty} \ n_j \ \wlr{\wfc{B_{n_j}}{f,x} - \wfc{f}{x}}$
converges to $L \in \overline{\wrone{}}$. According to Equation \pRef{diffN},
\[
\lim_{j \rightarrow \infty} \ n_j \ \wlr{\wfc{B_{n_j}}{f,x} - \wfc{f}{x}} =
\lim_{j \rightarrow \infty} \ \frac{n_j}{\wfc{D_{n_j}}{x}}
\wlr{\wfc{\Delta_{n_j}}{f,x} + \wfc{O}{f,x}} = L.
\]
Theorem \ref{thm2} in Section \ref{secNum} shows that
\[
\lim_{j \rightarrow \infty} \wfc{\Delta_{n_j}}{f,x} = 0,
\]
and if $\wfc{O}{f,x} = 0$ then $L = 0$, because
the sequence $n_j/\wfc{D_{n_j}}{x}$ is bounded by
Lemma \ref{lemDen} in Section \ref{secDen}. Since $0 \in \wfc{\wcal{O}}{f,x}$,
we are done when $\wfc{O}{f,x} = 0$ .
Let us then assume that $\wfc{O}{f,x} \neq 0$.
The equations above imply that
\[
\lim_{j \rightarrow \infty} \wabs{\frac{n_j}{\wfc{D_{n_j}}{x}}} = \wabs{L/\wfc{O}{f,x}},
\]
and Lemma \ref{lemDen} shows that
\[
\lim_{j \rightarrow \infty} \wabs{ \frac{n_j}{\wfc{D_{n_j}}{x}}} =
\lim_{j \rightarrow \infty}
\frac{
n_j \wfc{A}{\wfc{\rho_n^2}{x}}
}
{
\wabs{\wfc{D_{n_j}}{x}}
}
\frac{1}{\wfc{A}{\wfc{\rho_n^2}{x}}}
 =
\lim_{j \rightarrow \infty} \frac{1}{\wfc{A}{\wfc{\rho_n^2}{x}}} \leq 2 / \pi.
\]
Therefore, $\wabs{L} \leq 2 \wabs{\wfc{O}{f,x}}/\pi$ and Equation \pRef{ofIrrat} is correct.
Conversely, if $z \in \wfc{\wcal{O}}{f,x}$ then either $z = 0$ or
$\wabs{z} \in (0, 2 \wfc{O}{f,x}/\pi]$. In the
first case, Lemma \ref{lemIrrat} yields an odd sequence $n_j$ such that
\[
\lim_{j \rightarrow \infty} \frac{1}{n_j} \wfc{D_{n_j}}{x} = +\infty
\]
and we have that
\[
\lim_{j \rightarrow \infty} \ n_j \ \wlr{\wfc{B_{n_j}}{f,x} - \wfc{f}{x}} =
 \wfc{\Delta_n}{f,z} \frac{n_j}{\wfc{D_{n_j}}{x}}  =
 \wfc{O}{f,x} \times 0 = z.
\]
Otherwise, when $z \neq 0$,
\pbTClaim{ohMy}
y = z / \wfc{O}{f,x} \in [-2/\pi,2/\pi] \setminus \wset{0}
\peTClaim{ohMy}
and Lemma \ref{lemIrrat} yields an odd sequence $n_j$ such that
\[
\lim_{j \rightarrow \infty} \frac{1}{n_j} \wfc{D_{n_j}}{x} = 1/y,
\]
and
\[
\lim_{j \rightarrow \infty} \ n_j \ \wlr{\wfc{B_{n_j}}{f,x} - \wfc{f}{x}} =
 \wfc{\Delta_n}{f,z} \frac{n_j}{\wfc{D_{n_j}}{x}}  =
 \wfc{O}{f,x} \times z/ \wfc{O}{f,x} = z.
\]
Therefore, we have proved the converse part of Theorem \ref{thmMain} for
an irrational $x$ and an odd sequence $n_j$.
The same argument applies for an irrational $x$ and an even sequence $n_j$,
replacing $\wfc{O}{f,x}$ by $\wfc{E}{f,x}$
and $\wfc{\wcal{O}}{f,x}$ by $\wfc{\wcal{E}}{f,x}$.

Let us then analyze a rational $x$. Since $x$ is regular,
we must have $x \in \wlr{-1,1}$, and there exist
positive integers $p$ and $q$ with $\wfc{\wrm{gcd}}{p,q} = 1$
such that $x = p/q - 1$, and we can use
the argument applied in the irrational case
replacing the interval $[-2 /\pi, 2/\pi]$ in Equation
\pRef{ohMy} by the
set $\wfc{\wcal{O}}{p/q}$
or $\wfc{\wcal{E}}{p/q}$ in Corollaries \ref{corRatPQ},  \ref{corRat2PQ} and \ref{corRatP2Q}
in Section \ref{secDen} corresponding to the parity of $p$ and $q$,
and replacing the intervals
\[
\left[\ -2 \  \wabs{\wfc{O}{f,x}}/\pi, 2 \ \wabs{\wfc{O}{f,x}}/\pi \ \right]
\hspace{1cm} \wrm{and} \hspace{1cm}
\left[ \ -2 \ \wabs{\wfc{E}{f,x}}/\pi, 2 \ \wabs{\wfc{E}{f,x}}/\pi \ \right]
\]
by the sets
\pbDef{oddSetRat}
\wfc{\wcal{O}}{f,x} =
\wset{\wfc{O}{f,x} / y, \ \ y \in \wfc{\wcal{O}}{p/q}}
\peDef{oddSetRat}
and
\pbDef{evenSetRat}
\wfc{\wcal{E}}{f,x} = \wset{ \wfc{E}{f,x}/y, \ \ y \in \wfc{\wcal{E}}{p/q}}.
\peDef{evenSetRat}
\peProof{Theorem}{thmMain}\\

\pbProofB{Theorem}{thmUnif}
Theorem \ref{thmUnif} follows from Lemma \ref{lemDen}
and Theorem \ref{thm2}.
\peProof{Theorem}{thmUnif}\\

\pbProofB{Theorem}{thmUnifBV}
Equations \pRef{diffNO} and \pRef{diffNE} show that
\[
n \wnorm{\wfc{B_n}{f} - f}_{\infty}
\leq \frac{
\wnorm{\wfc{\Delta_n}{f}}_{\infty} + 
\max \wset{\wnorm{\wfc{O}{f}}_{\infty}, \wnorm{\wfc{E}{f}}_{\infty}}} 
{
\wnorm{D_n}_{\infty} / n
},
\]
and Equation \pRef{boundDenA} and Corollary \ref{corNumBV} imply
Equation \pRef{unifBV}.
\peProof{Theorem}{thmUnifBV}\\

\end{document}